# FEEDBACK STABILIZATION METHODS FOR THE NUMERICAL SOLUTION OF SYSTEMS OF ORDINARY DIFFERENTIAL EQUATIONS

**Iasson Karafyllis[*] and Lars Grüne[**]**

[*]**Department of Environmental Engineering,**
Technical University of Crete,73100, Chania, Greece
email: ikarafyl@enveng.tuc.gr

[**]**Mathematisches Institut, Fakultät für Mathematik und Physik,**
Universität Bayreuth, 95440 Bayreuth, Germany
email: lars.gruene@uni-bayreuth.de

**Abstract**
In this work we study the problem of step size selection for numerical schemes, which guarantees that the numerical solution presents the same qualitative behavior as the original system of ordinary differential equations, by means of tools from nonlinear control theory. Lyapunov-based and Small-Gain feedback stabilization methods are exploited and numerous illustrating applications are presented for systems with a globally asymptotically stable equilibrium point. The obtained results can be used for the control of the global discretization error as well.

## 1. Introduction

It is well-known that step size control can enhance the performance of a numerical scheme when applied to a system of Ordinary Differential Equations (ODEs). In fact the use of the word "control" suggests that methods and techniques used in Mathematical Control Theory can be (in principle) used in order to achieve certain objectives for the numerical solution of systems of ODEs. For example, in [6] the authors use a "Proportional-Integral" technique which is similar to the "Proportional-Integral" controller used in Linear Systems Theory in order to keep the local discretization error within certain bounds (see also [7,8,13]). Theoretical results on the behavior of adaptive time-stepping methods have been presented in [20,22] and the control theoretic notion of input-to-state stability (ISS) has been successfully used in [10,11] in order to explain the behavior of attractors under discretization.

In this work, we develop tools for nonlinear systems which are similar to methods used in Nonlinear Control Theory. We consider the problem of selecting the step size for numerical schemes so that the numerical solution presents the same qualitative behavior as the original system of ODEs. It is well-known that any consistent and stable numerical scheme for ODEs inherits the asymptotic stability of the original equation in a practical sense, even for more general attractors than equilibria see for instance [10,11] and [28] (Chapter 7). Practical Asymptotic stability means that the system exhibits an asymptotically stable set close to the original attractor, i.e., in our case a small neighbourhood around the equilibrium point, which shrinks down to the attractor as the time step $h$ tends to 0. In contrast to these results, in this paper we investigate the case in which the numerical approximation is asymptotically stable in the usual sense, i.e., not only practically.

Here, we concentrate on nonlinear systems for which an equilibrium point is the global attractor. In Section 2 of the present work, it is shown how the problem of appropriate step size selection can be converted to a rigorous abstract feedback stabilization problem for a particular hybrid system (see also [17]-the reader should notice that the standard stability analysis of numerical schemes uses the study of a discrete-time system e.g., [13,14,16,21,28]; not a hybrid system). Therefore, we are in a position to use all methods of feedback design for nonlinear systems. Specifically, we consider

- ➢ methods based on Small-Gain Theorems,
- ➢ methods based on Lyapunov functions.

Both methods have been used widely in Nonlinear Systems Theory for the solution of feedback stabilization problems (see [1,3,15,18,19,25,27] and references therein). In the present work, the above methods are used for the step size selection for numerical schemes for ODEs (see Section 3 and Section 4). General results are developed for arbitrary consistent Runge-Kutta schemes (see Theorem 4.5, Theorem 4.9 and Theorem 4.12 below) and specific results are given for specific Runge-Kutta schemes (see Corollary 4.7 and Theorem 4.15 below). The analysis is based on the Control Lyapunov Function methodology. The obtained results constitute nonlinear extensions of well-known properties of numerical schemes (e.g., A-stability, cf. Corollary 4.16). On the other hand small-gain methods allow the proposal of novel numerical schemes (see Theorem 3.1 below) for certain classes of nonlinear systems.

A number of applications of the obtained results is developed in Sections 5, 6 and 7. More specifically, we consider the possibility of using explicit schemes for stiff linear systems of ODEs: our results constitute rigorous nonlinear extensions of the ideas presented in [2] (see Section 6). Furthermore, we consider the problem of controlling the **global** discretization error: our results show that rigorous global discretization error control is possible **only after** stabilization of the numerical scheme (by appropriate step size selection-see Section 7). Another application of stabilization methods based on Small-Gain analysis to systems described by Partial Differential Equations is presented in Section 5.

Thus, the contribution of the paper is twofold. On the one hand, our control theoretic approach yields new insight into the stability properties of numerical schemes and as such it adds another means to the toolbox for stability investigations of numerical schemes. On the other hand, our method leads to the design of new discretization schemes and step size control algorithms, which beyond the mere control of the local discretization error also take care of the global qualitative behaviour.

**Notations** Throughout this paper we adopt the following notations:

∗ Let $A \subseteq \Re^n$ be a set. By $C^0(A;\Omega)$, we denote the class of continuous functions on $A$, which take values in $\Omega$. By $C^k(A;\Omega)$, where $k \geq 1$ is an integer, we denote the class of differentiable functions on $A$ with continuous derivatives up to order $k$, which take values in $\Omega$. By $C^\infty(A;\Omega)$, we denote the class of differentiable functions on $A$ having continuous derivatives of all orders, which take values in $\Omega$, i.e., $C^\infty(A;\Omega) = \bigcap_{k \geq 1} C^k(A;\Omega)$.

∗ For a vector $x \in \Re^n$ we denote by $|x|$ its usual Euclidean norm and by $x'$ its transpose. By $B_\varepsilon(x)$, where $\varepsilon > 0$ and $x \in \Re^n$, we denote the ball of radius $\varepsilon > 0$ centered at $x \in \Re^n$, i.e., $B_\varepsilon(x) := \{ y \in \Re^n : |y - x| < \varepsilon \}$. For a real matrix $A \in \Re^{n \times m}$ we denote by $|A|$ its induced norm, i.e., $|A| := \max\{ |Ax| : x \in \Re^m, |x| = 1 \}$ and by $A' \in \Re^{m \times n}$ its transpose.

∗ $\Re^+$ denotes the set of non-negative real numbers and $Z^+$ the set of non-negative integer numbers. $C$ denotes the set of complex numbers.

∗ By $K_\infty$ we denote the set of all increasing and continuous functions $\rho : \Re^+ \to \Re^+$ with $\rho(0) = 0$ and $\lim_{s \to +\infty} \rho(s) = +\infty$. By $KL$ we denote the set of all continuous functions $\sigma = \sigma(s,t) : \Re^+ \times \Re^+ \to \Re^+$ with the properties: (i) for each $t \geq 0$ the mapping $\sigma(\cdot,t)$ is of class $K$; (ii) for each $s \geq 0$, the mapping $\sigma(s,\cdot)$ is non-increasing with $\lim_{t \to +\infty} \sigma(s,t) = 0$.

∗ For every scalar continuously differentiable function $V : \Re^n \to \Re$, $\nabla V(x)$ denotes the gradient of $V$ at $x \in \Re^n$, i.e., $\nabla V(x) = \left( \frac{\partial V}{\partial x_1}(x), \ldots, \frac{\partial V}{\partial x_n}(x) \right)$. We say that a function $V : \Re^n \to \Re^+$ is positive definite if $V(x) > 0$ for all $x \neq 0$ and $V(0) = 0$. We say that a continuous function $V : \Re^n \to \Re^+$ is radially unbounded if the following property holds: "for every $M > 0$ the set $\{ x \in \Re^n : V(x) \leq M \}$ is compact".

∗ For a sufficiently smooth function $V : \Re^n \to \Re$ we denote by $L_f V(x) := \nabla V(x) f(x)$ the Lie derivative of $V$ along $f$ and we define recursively $L_f^{(i+1)} V(x) = L_f \left( L_f^{(i)} V(x) \right)$ for $i \geq 1$.

## 2. Description of the Problem

Consider the autonomous system

$$\dot{z}(t) = f(z(t)), \; z(t) \in \Re^n \tag{2.1}$$

where $f : \Re^n \to \Re^n$ is a locally Lipschitz vector field with $f(0) = 0$. For every $z_0 \in \Re^n$ and $t \geq 0$, the solution of (2.1) with initial condition $z(0) = z_0$ will be denoted by $z(t, z_0)$. The numerical approximation of system (2.1) will be the hybrid system:

$$\begin{aligned}
&\dot{x}(t) = F(h_i, x(\tau_i)), \; t \in [\tau_i, \tau_{i+1}) \\
&\tau_0 = 0, \; \tau_{i+1} = \tau_i + h_i \\
&h_i = \varphi(x(\tau_i)) \exp(-u(\tau_i)) \\
&x(t) \in \Re^n, \; u(t) \in [0, +\infty)
\end{aligned} \tag{2.2}$$

where $\varphi \in C^0(\Re^n; (0, r])$, $r > 0$ is a constant, $F : \bigcup_{x \in \Re^n} ([0, \varphi(x)] \times \{x\}) \to \Re^n$ is a (not necessarily continuous) vector field with $F(h, 0) = 0$ for all $h \in [0, \varphi(0)]$, $\lim_{h \to 0^+} F(h, z) = f(z)$, for all $z \in \Re^n$. More specifically, the solution $x(t)$ of the hybrid system (2.2) is obtained for every locally bounded $u : \Re^+ \to \Re^+$ and $x_0 \in \Re^n$ by the following algorithm (see [17]):

Step $i$ :
1) Given $\tau_i$ and $x(\tau_i)$, calculate $\tau_{i+1}$ using the equation $\tau_{i+1} = \tau_i + \varphi(x(\tau_i)) \exp(-u(\tau_i))$,
2) Compute the state trajectory $x(t)$, $t \in (\tau_i, \tau_{i+1}]$ as the solution of the differential equation $\dot{x}(t) = F(h_i, x(\tau_i))$, i.e., $x(t) = x(\tau_i) + (t - \tau_i) F(h_i, x(\tau_i))$ for $t \in (\tau_i, \tau_{i+1}]$.

For $i = 0$ we take $\tau_0 = 0$ and $x(0) = x_0$ (initial condition).

We will further assume that there exists a continuous, non-decreasing function $M : \Re^+ \to \Re^+$ such that

$$|F(h, x)| \leq |x| M(|x|) \text{ for all } x \in \Re^n \text{ and } h \in [0, \varphi(x)] \tag{2.3}$$

It should be noticed that the hybrid system (2.2) under hypothesis (2.3) is an autonomous system, which satisfies the "Boundedness-Implies-Continuation" property and for each locally bounded input $u : \Re^+ \to \Re^+$ and $x_0 \in \Re^n$ there exists a unique absolutely continuous function $[0, +\infty) \ni t \to x(t) \in \Re^n$ with $x(0) = x_0$ which satisfies (2.2) (see [17]). Some remarks are needed in order to justify the name "numerical approximation of system (2.1)" for the hybrid system (2.2):

a) Notice that the condition $\lim_{h \to 0^+} F(h, z) = f(z)$ is a consistency condition for the numerical scheme applied to (2.1).
b) The sequence $\{h_i\}_0^\infty$ is the sequence of step sizes used in order to obtain the numerical solution. Notice that for the case $\varphi(x) \equiv r$, constant inputs $u(t) \equiv u \geq 0$ will produce constant step sizes with $h_i \equiv r \exp(-u)$. Moreover, notice that variable step sizes can be represented easily, by selecting in an appropriate way the input $u : \Re^+ \to \Re^+$.
c) The constant $r > 0$ is the maximum allowable step size.
d) The function $\varphi \in C^0(\Re^n; (0, r])$ determines the maximum allowable step size $\varphi(x(\tau_i))$ for each $x(\tau_i) \in \Re^n$. This is important for implicit numerical schemes as shown below.

All consistent $s$–stage Runge-Kutta methods can be represented by the hybrid system (2.2). More specifically, let $x_0 \in \Re^n$ and consider a consistent $s$–stage Runge-Kutta method for (2.1):

$$Y_i = x_0 + h\sum_{j=1}^{s} a_{ij} f(Y_j) \quad , \quad i = 1,...,s \tag{2.4}$$

$$x = x_0 + h\sum_{i=1}^{s} b_i f(Y_i) \tag{2.5}$$

with $\sum_{i=1}^{s} b_i = 1$. If the scheme is explicit, i.e., if $a_{ij} = 0$ for $j \geq i$, then there always exists a unique solution to equations (2.4). If the scheme is implicit, then in order to be able to guarantee that equations (2.4) admit a unique solution it may be necessary to restrict the step size to $h \in [0, \varphi(x_0)]$ for some maximal step size $\varphi(x_0)$ depending on the state $x_0 \in \Re^n$. In all subsequent statements on numerical schemes, we will tacitly assume that such a step size restriction is imposed if necessary.

A suitable choice for $\varphi(x)$ may be obtained in the following way. Let $\gamma : \Re^+ \to \Re^+$ be a continuous, non-decreasing function with $|f(x)| \leq |x|\gamma(|x|)$ for all $x \in \Re^n$ (such a function always exists since $f : \Re^n \to \Re^n$ is a locally Lipschitz vector field with $f(0) = 0$). Let $L_\lambda : \Re^n \to (0,+\infty)$ be a continuous function with $L_\lambda(x_0) \geq \sup\left\{\frac{|f(x) - f(y)|}{|x - y|} : x, y \in B_\lambda(x_0), x \neq y\right\}$ for all $x_0 \in (\Re^n \setminus \{0\})$, with $B_\lambda(x_0) := \{x \in \Re^n : |x - x_0| \leq \lambda|x_0|\}$, $\lambda \in (0,1)$. The continuous function $\varphi(x) := \frac{\lambda}{|A|(L_\lambda(x) + \gamma(|x|))}$, where $|A| := \max_{i=1,...,s} \sum_{j}^{s} |a_{ij}|$, guarantees that for all $x_0 \in \Re^n$ and $h \in [0, \varphi(x_0)]$ the equations (2.4) have a unique solution satisfying $Y_i \in B_\lambda(x_0)$, $i = 1,...,s$. Note however that this bound may be conservative. For instance, if we apply the implicit Euler scheme ($s = 1, a_{11} = b_1 = 1$) to an asymptotically stable linear ODE of the form $\dot{x} = J x$ with a Hurwitz matrix $J \in \Re^{n \times n}$, then (2.4) becomes

$$Y_1 = x_0 + hJY_1 \Leftrightarrow (I - hJ)Y_1 = x_0$$

which always has a unique solution because all eigenvalues of $-J$ and thus of $I - hJ$ have positive real parts for all $h \geq 0$; hence $I - hJ$ is invertible for all $h \geq 0$.

We define

$$F(h, x_0) := h^{-1}(x - x_0) = \sum_{i=1}^{s} b_i f(Y_i) \tag{2.6}$$

A moment's thought reveals that for every locally bounded $u : \Re^+ \to \Re^+$ and $x_0 \in \Re^n$ the solution of (2.2) with (2.6) coincides at each $\tau_i$, $i \geq 0$ with the numerical solution of (2.1) with $x(0) = x_0$ obtained by using the Runge-Kutta numerical scheme (2.4), (2.5) and using the discretization step sizes $h_i = \varphi(x(\tau_i))\exp(-u(\tau_i))$, $i \geq 0$. The reader should notice that other ways (besides (2.6)) of defining the vector field $F : \bigcup_{x \in \Re^n} ([0, \varphi(x)] \times \{x\}) \to \Re^n$ may be possible: here we have selected the simplest way of obtaining a piecewise linear numerical solution.

Moreover, the reader should notice that appropriate step size restriction can always guarantee that (2.3) holds for $F : \bigcup_{x \in \Re^n} ([0, \varphi(x)] \times \{x\}) \to \Re^n$ as defined by (2.6). For example, if $\varphi(x) := \frac{\lambda}{|A|(L_\lambda(x) + \gamma(|x|))}$ is the step size restriction described above, then $F : \bigcup_{x \in \Re^n} ([0, \varphi(x)] \times \{x\}) \to \Re^n$ as defined by (2.6) satisfies $|F(h,x)| \leq |x|\left[1 + r(1+\lambda)\left(\sum_{i=1}^{s} |b_i|\right)\gamma((1+\lambda)|x|)\right]$ for all $x \in \Re^n$ and $h \in [0, \varphi(x)]$. Thus (2.3) holds with

$$M(y) := 1 + r(1+\lambda)\left(\sum_{i=1}^{s}|b_i|\right)\gamma((1+\lambda)y).$$

If the Runge-Kutta scheme (2.4), (2.5) is of order $p \geq 1$, we will occasionally further assume that $f \in C^p(\Re^n; \Re^n)$ and for each fixed $x \in \Re^n$ the mapping $[0, \varphi(x)] \ni h \to F(h, x)$ is $p$ times continuously differentiable with

$$|F(h,x)| + \sum_{j=1}^{p}\left|\frac{\partial^j}{\partial h^j}F(h,x)\right| \leq G(|x|)\max\left\{|f(y)| : y \in \Re^n, |y-x| \leq |x|\varphi(x)M(|x|)\right\} \text{ for all } x \in \Re^n \text{ and } h \in [0, \varphi(x)] \quad (2.7)$$

for some continuous, non-decreasing function $G: \Re^+ \to \Re^+$, where $M: \Re^+ \to \Re^+$ is the function involved in (2.3). The reader should notice that appropriate step size restriction can always guarantee that (2.7) holds for $F: \bigcup_{x \in \Re^n}([0,\varphi(x)] \times \{x\}) \to \Re^n$ as defined by (2.6). Notice that the implicit function theorem for (2.4) guarantees for each fixed $x \in \Re^n$ the existence of $\varphi(x) > 0$ such that the mapping $[0, \varphi(x)] \ni h \to F(h,x)$ is $p$ times continuously differentiable. A suitable choice for $\varphi(x)$ may be obtained by the formula $\varphi(x) := \frac{\lambda}{1 + 2|A|\max\{|Df(z)| : |z| \leq (1+\lambda)|x|\}}$, where $\lambda \in (0,1)$, $|A| := \max_{i=1,\ldots,s}\sum_{j}|a_{ij}|$. However, it must be emphasized again that the above step size restriction may be conservative (e.g., for explicit schemes).

Using Theorem II.3.1 in [12], (2.7), the fact that $f \in C^p(\Re^n; \Re^n)$ and the fact that $g_k(z(h,x)) = \frac{\partial^k}{\partial h^k}z(h,x)$ for $k \geq 1$, where $g_k : \Re^n \to \Re^n$ for $k = 1, \ldots, p+1$ are vector fields obtained by the recursive formulae $g_1(z) = f(z)$, $g_{i+1}(z) = Dg_i(z)f(z)$, we may conclude that there exist continuous functions $N : \Re^n \to (0, +\infty)$, $C : \Re^n \to \Re^+$ such that the following inequalities hold for all $x \in \Re^n$ and $h \in [0, \varphi(x)]$:

$$C(x) \leq N(x)\left[\max\{|f(y)| : y \in \Re^n, |y-x| \leq |x|\varphi(x)M(|x|)\} + \max\{|f(z(h,x))| : h \in [0, \varphi(x)]\}\right] \quad (2.8a)$$

$$|z(h,x) - x - hF(h,x)| \leq h^{p+1}C(x) \quad (2.8b)$$

If we further assume that there exists a neighborhood $N \subseteq \Re^n$ with $0 \in N$ satisfying the following properties:

➢ there exists a constant $\Lambda > 0$ and an integer $q \geq 1$ such that $|f(x)| \leq \Lambda|x|^q$ for all $x \in N$,
➢ there exists a constant $Q > 0$ such that $|z(h,x)| \leq Q|x|$ for all $x \in N$ and $h \in [0, \varphi(x)]$,

then it follows from (2.8a) that there exists a neighborhood $\tilde{N} \subseteq N$ with $0 \in \tilde{N}$ and a constant $K > 0$ such that

$$C(x) \leq Kh^{p+1}|x|^q, \text{ for all } x \in \tilde{N} \quad (2.9)$$

Assume next that $0 \in \Re^n$ is Uniformly Globally Asymptotically Stable (UGAS) for (2.1) (in the sense described in [24]; see also [17,19]). Our goal is to be able to produce numerical solutions using the numerical approximation (2.2) which have the correct qualitative behavior. More specifically, we would like to be in a position to know a continuous function $\varphi : \Re^n \to (0, r]$ so that the numerical solution produced by (2.2) has the correct qualitative behavior (e.g., $\lim_{t \to +\infty} x(t) = 0$). However, we would like to be able to guarantee that the correct behavior for the numerical solution can be obtained by using arbitrary discretization step sizes smaller than $\varphi(x(\tau_i))$ (i.e., if we obtain the correct qualitative behavior using the discretization step sizes $h_i = \varphi(x(\tau_i))$ $i \geq 0$, we would like to obtain the correct qualitative behavior using the discretization step sizes $h_i = \varphi(x(\tau_i))\exp(-u(\tau_i))$, $i \geq 0$, where $u : \Re^+ \to \Re^+$ is an arbitrary locally bounded function). This is equivalent by requiring that $0 \in \Re^n$ is Uniformly

Robustly Globally Asymptotically Stable (URGAS) for (2.2) (in the sense described in [17]).

The reader should notice that continuity for the function $\varphi: \Re^n \to (0, r]$ is essential: without assuming continuity it may happen that $\liminf_{x \to 0} \varphi(x) = 0$ and this would require discretization step sizes of vanishing magnitude as $t \to +\infty$.

Moreover, since we want to be able to determine a continuous function $\varphi: \Re^n \to (0, r]$, which "stabilizes" the hybrid system (2.2), we are essentially studying a feedback stabilization problem for the hybrid system (2.2). Hence, we are in a position to pose the problem rigorously. We consider the following feedback stabilization problems:

**(P1)-Existence Problem:** *Is there a continuous function $\varphi: \Re^n \to (0, r]$, such that $0 \in \Re^n$ is URGAS for system (2.2)?*

**(P2)-Design Problem:** *Construct a continuous function $\varphi: \Re^n \to (0, r]$, such that $0 \in \Re^n$ is URGAS for system (2.2).*

It is well known that any consistent and stable numerical scheme for ODEs inherits the asymptotic stability of the original equation in a practical sense, even for more general attractors than equilibria see for instance [10,11] or [28] (Chapter 7). Practical Asymptotic stability means that the system exhibits an asymptotically stable set close to the original attractor, i.e., in our case a small neighbourhood around the equilibrium point, which shrinks down to the attractor as the time step $h$ tends to 0.

Here, the property we are looking for, i.e., "real" asymptotic stability, is a stronger property which cannot in general be deduced from practical stability. In [28] (Chapter 5), several results for our problem for specific classes of ODEs are derived using classical numerical stability concepts like A-stability, B-stability and the like. In contrast to this reference, in the sequel we use nonlinear control theoretic analysis and feedback design techniques; more precisely Small-Gain and Lyapunov function techniques in Sections 3 and 4, respectively for solving Problems (P1) and (P2). This allows us to obtain asymptotic stability results under different structural assumptions and for more general classes of systems as in [28] (Chapter 5).

## 3. Small-Gain Methodology

One of the tools used in mathematical control theory for nonlinear feedback design is the methodology based on small-gain results. The methodology was first used in [15] where a nonlinear small-gain result based on the notion of Input-to-State Stability (ISS-see [26]) was presented. Since then it has been applied successfully to many feedback stabilization problems. Recently, the small-gain theorem was extended to general control systems including hybrid systems (see [18]). Consequently, the small-gain methodology for feedback design appears to be applicable for the solution of problem (P2) for certain classes of nonlinear systems (2.1).

Consider the following system:

$$\dot{z} = f_0(z) \quad , \quad z \in \Re^m \tag{3.1a}$$

$$\begin{aligned} \dot{x}_1 &= -a_1(x_1)x_1 + f_1(z) \\ \dot{x}_i &= -a_i(x_i)x_i + f_i(z, x_1, \ldots, x_{i-1}) \quad , \quad i = 2, \ldots, n \\ x &= (x_1, \ldots, x_n)' \in \Re^n \end{aligned} \tag{3.1b}$$

where $f_0: \Re^m \to \Re^m$, $f_1: \Re^m \to \Re$, $f_i: \Re^m \times \Re^{i-1} \to \Re$ ($i = 2, \ldots, n$), $a_i: \Re \to \Re$ ($i = 1, \ldots, n$) are locally Lipschitz mappings with $f_0(0) = 0$, $f_1(0) = \ldots = f_n(0, 0, \ldots, 0) = 0$. We assume that there exist constants $L_i > 0$ ($i = 1, \ldots, n$) such that:

$$a_i(y) \geq L_i, \quad \forall y \in \Re \tag{3.2}$$

We also assume that $0 \in \Re^m$ is UGAS for (3.1a). Under the previous assumptions, using the fact that system (3.1) has a structure of systems in cascade, we may prove by induction that for every $j = 1, \ldots, n$, $0 \in \Re^m \times \Re^j$ is UGAS for

system (3.1a) with

$$\begin{aligned}\dot{x}_1 &= -a_1(x_1)x_1 + f_1(z) \\ \dot{x}_i &= -a_i(x_i)x_i + f_i(z, x_1,..., x_{i-1}) \quad , \quad i = 2,..., j\end{aligned} \tag{3.3}$$

The proof for $j=1$ is based on the fact that for every $x_{10} \in \Re$ and for every measurable $u : \Re^+ \to \Re$ the solution of $\dot{x}_1 = -a_1(x_1)x_1 + u$ with initial condition $x_1(0) = x_{10}$ satisfies the following estimate:

$$|x_1(t)| \le \exp\left(-\frac{L_1}{2}t\right)|x_{10}| + \frac{1}{L_1}\sup_{0 \le s \le t}|u(s)|, \quad \forall t \ge 0 \tag{3.4}$$

Consequently, the solution of $\dot{x}_1 = -a_1(x_1)x_1 + f_1(z)$ satisfies $|x_1(t)| \le \exp\left(-\frac{L_1}{2}t\right)|x_{10}| + \frac{1}{L_1}\sup_{0 \le s \le t}|f_1(z(s))|$, i.e., satisfies the uniform ISS property from the input $z \in \Re^m$. Since $0 \in \Re^m$ is GAS for (3.1a), a well-known corollary of the small-gain theorem (systems in cascade) guarantees UGAS for the composite system. The proof is similar for all $j = 1,..., n$.

Suppose that a stable numerical scheme is available for (3.1a), i.e., there exist $\varphi \in C^0(\Re^m; (0, r])$, $r > 0$ and $F_0 : \bigcup_{z \in \Re^m}([0, \varphi(z)] \times \{z\}) \to \Re^m$ with $F_0(h, 0) = 0$ for all $h \in [0, \varphi(0)]$, $\lim_{h \to 0^+} F_0(h, z) = f_0(z)$, for all $z \in \Re^m$ such that $0 \in \Re^m$ is URGAS for the hybrid system:

$$\begin{aligned}\dot{z}(t) &= F_0(h_i, z(\tau_i)), \, t \in [\tau_i, \tau_{i+1}) \\ \tau_0 &= 0, \, \tau_{i+1} = \tau_i + h_i \\ h_i &= \varphi(z(\tau_i))\exp(-u(\tau_i)) \\ x(t) &\in \Re^n, \, u(t) \in [0, +\infty)\end{aligned} \tag{3.5}$$

We propose the following first order numerical scheme for (3.1b):

$$\begin{aligned}x_1(t+h) &= x_1(t) - ha_1(x_1(t))x_1(t+h) + hf_1(z(t)) \\ x_i(t+h) &= x_i(t) - ha_i(x_i(t))x_i(t+h) + hf_i(z(t), x_1(t),..., x_{i-1}(t)) \quad , \quad i = 2,..., n\end{aligned} \tag{3.6}$$

The above scheme is a partitioned scheme which treats differently the state $x_i$ and the states $z, x_1,..., x_{i-1}$ for each differential equation. The resulting hybrid system is system (3.5) with:

$$\begin{aligned}\dot{x}_1(t) &= \frac{-a_1(x_1(\tau_i))}{1 + h_i a_1(x_1(\tau_i))} x_1(\tau_i) + \frac{1}{1 + h_i a_1(x_1(\tau_i))} f_1(z(\tau_i)) \\ \dot{x}_j(t) &= \frac{-a_j(x_j(\tau_i))}{1 + h_i a_j(x_j(\tau_i))} x_j(\tau_i) + \frac{1}{1 + h_i a_j(x_j(\tau_i))} f_j(z(\tau_i), x_1(\tau_i),..., x_{j-1}(\tau_i)) \quad , \quad j = 2,..., n\end{aligned} \tag{3.7}$$

We have:

**Theorem 3.1:** $0 \in \Re^m \times \Re^n$ *is uniformly RGAS for system (3.5), (3.7).*

The proof of the above theorem is based on the following technical lemma.

**Lemma 3.2:** *Let $a : \Re \to \Re$ be a continuous function with $L = \inf_{y \in \Re} a(y) > 0$ and let a constant $r > 0$. Then for every sequence $\{h_i\}_0^\infty$ with $h_i \in (0, r]$ for all $i \ge 0$, for every locally bounded function $v : \Re^+ \to \Re$ and for every $x_0 \in \Re$ the solution of*

$$\dot{x}(t) = \frac{-a(x(\tau_i))}{1+h_i a(x(\tau_i))} x(\tau_i) + \frac{1}{1+h_i a(x(\tau_i))} v(\tau_i) \quad , \quad t \in [\tau_i, \tau_{i+1}) \tag{3.8}$$

$$\tau_{i+1} = \tau_i + h_i \, , h_i \in (0, r], x(t) \in \Re$$

with initial condition $x(0) = x_0 \in \Re$, $\tau_0 = 0$ satisfies the following estimate

$$|x(t)| \leq \exp(\sigma r)|x_0|\exp(-\sigma t) + \frac{1+e}{e\sigma L} \sup_{0 \leq s \leq t} |v(s)|, \quad \forall t \in \left[0, \sup_{i \geq 0} \tau_i \right) \tag{3.9}$$

where $\sigma > 0$ is any constant such that $\dfrac{1}{1+s} \leq \exp(-\sigma s)$ for all $s \in [0, rL]$.

**Proof:** Notice that for every $i \geq 0$ it holds that:

$$x(\tau_{i+1}) = x_0 \prod_{j=0}^{i}(1+h_j a(x(\tau_j)))^{-1} + \sum_{j=0}^{i}\left[h_j v(\tau_j)\left(\prod_{k=j}^{i}(1+h_k a(x(\tau_k)))^{-1}\right)\right] \tag{3.10}$$

and using the definition $L = \inf_{y \in \Re} a(y) > 0$, we obtain the following bound from (3.10):

$$|x(\tau_{i+1})| \leq |x_0|\prod_{j=0}^{i}(1+h_j L)^{-1} + \max_{j=0,\ldots,i}|v(\tau_j)|\sum_{j=0}^{i}\left[h_j\left(\prod_{k=j}^{i}(1+h_k L)^{-1}\right)\right] \tag{3.11}$$

Let $\sigma > 0$ such that $\dfrac{1}{1+s} \leq \exp(-\sigma s)$ for all $s \in [0, rL]$. It follows that

$$\prod_{j=0}^{i}(1+h_j L)^{-1} \leq \prod_{j=0}^{i}\exp(-\sigma L h_j) = \exp(-\sigma L \tau_{i+1})$$

and

$$\sum_{j=0}^{i}\left[h_j\left(\prod_{k=j}^{i}(1+h_k L)^{-1}\right)\right] \leq \sum_{j=0}^{i}\left[h_j\left(\prod_{k=j}^{i}\exp(-\sigma L h_k)\right)\right] = \sum_{j=0}^{i}\left[h_j \exp(-\sigma L(\tau_{i+1}-\tau_j))\right] =$$

$$= \exp(-\sigma L \tau_{i+1})\sum_{j=0}^{i}\left[\exp(\sigma L \tau_j)\int_{\tau_j}^{\tau_{j+1}} ds\right] \leq \exp(-\sigma L \tau_{i+1})\sum_{j=0}^{i}\left[\int_{\tau_j}^{\tau_{j+1}}\exp(\sigma L s)ds\right] =$$

$$= \exp(-\sigma L \tau_{i+1})\int_{0}^{\tau_{i+1}}\exp(\sigma L s)ds \leq \frac{1}{\sigma L}$$

Using the above inequalities in conjunction with (3.11) we obtain for all $i \geq 0$:

$$|x(\tau_{i+1})| \leq |x_0|\exp(-\sigma \tau_{i+1}) + \frac{1}{\sigma L}\max_{0 \leq j \leq i}|v(\tau_j)| \tag{3.12}$$

Notice that for every $i \geq 0$ and $t \in [\tau_i, \tau_{i+1})$ it holds that:

$$|x(t)| \leq |x(\tau_i)| + \frac{h_i}{1+h_i L}|v(\tau_i)| \tag{3.13}$$

Since $h_i \in (0, r]$, we obtain $\frac{1}{1+h_i L} \leq \exp(-\sigma h_i L)$ and consequently $\frac{h_i}{1+h_i L} \leq h_i \exp(-\sigma h_i L) \leq \frac{1}{e\sigma L}$. Combining (3.12) and (3.13) we obtain (3.13). The proof is complete. ◁

The proof of the Theorem 3.1 follows induction: we show that for every $k = 1,...,n$, $0 \in \Re^m \times \Re^k$ is uniformly RGAS for system (3.5) with

$$\dot{x}_1(t) = \frac{-a_1(x_1(\tau_i))}{1+h_i a_1(x_1(\tau_i))} x_1(\tau_i) + \frac{1}{1+h_i a_1(x_1(\tau_i))} f_1(z(\tau_i))$$

$$\dot{x}_j(t) = \frac{-a_j(x_j(\tau_i))}{1+h_i a_j(x_j(\tau_i))} x_j(\tau_i) + \frac{1}{1+h_i a_j(x_j(\tau_i))} f_j(z(\tau_i), x_1(\tau_i),...,x_{j-1}(\tau_i)) \quad , \quad j = 2,...,k$$

Notice that the Lemma 3.2 guarantees that $|x_i(t)| \leq \exp(\sigma r)|x_i(0)| \exp(-\sigma t) + \frac{1+e}{e\sigma L_i} \sup_{0 \leq s \leq t} |f_i(z(s),...,x_{i-1}(s))|$, where $\sigma > 0$ is a constant with the property $\frac{1}{1+s} \leq \exp(-\sigma s)$ for all $s \in [0, r \max_{i=1,...,n}(L_i)]$. Remark 3.2(b) in [18] (systems in cascade) guarantees URGAS for the composite system.

## 4. Lyapunov function based Step Selection

While the small-gain methodology can be applied to certain systems of differential equations and can yield numerical methods which guarantee numerical stability, it cannot be applied in a systematic way. On the other hand the Lyapunov-based feedback design methods can be applied to general nonlinear systems of differential equations and yield explicit formulas for the feedback law (see [25]). In this section we apply the Lyapunov-based feedback design methodology for the solution of Problems (P1) and (P2). It is well known that Lyapunov functions exist for every asymptotically stable ODE system and in many applications one can even give explicit formulas for these functions (some examples can be found in Section 6). However, even if a Lyapunov function is not exactly known, under suitable assumptions on the ODE system, certain structural properties of the Lyapunov function can be obtained (cf., e.g., Proposition 4.4, below) and used in our context. Hence, the main task of this section is to derive conditions under which the Lyapunov function for the ODE system can be used in order to conclude stability for the hybrid system (2.2), i.e., for the numerical approximation of system (2.1).

The results will be developed in the following way: first (subsection 4.I) we provide some background material needed for the derivation of the main results. In Subsection 4.II we consider general consistent Runge-Kutta schemes and provide sufficient conditions for the solvability of Problem (P1) and (P2). The results are specialized for the explicit Euler method. Finally, in Subsection 4.III, we present special results for the implicit Euler numerical scheme.

### 4.I. Background Material

The crucial technical result that allows the use of Lyapunov functions for hybrid systems of the form (2.2) is the following lemma. Here it should be recalled that $0 \in \Re^n$ is robustly K-exponentially stable for (2.2) if there exist a function $a \in K_\infty$ and a constant $\sigma > 0$ such that for every locally bounded $u : \Re^+ \to \Re^+$ and $x_0 \in \Re^n$ the solution $x(t, x_0; u)$ of (2.2) with initial condition $x(0) = x_0$ corresponding to $u : \Re^+ \to \Re^+$ satisfies the inequality $|x(t, x_0; u)| \leq \exp(-\sigma t) a(|x_0|)$ for all $t \geq 0$ (an extension of the corresponding notion for systems described by ODEs, see [23]).

**Lemma 4.1:** *Consider system (2.2) and suppose that there exist a continuous, positive definite and radially unbounded function $V : \Re^n \to \Re^+$ and a continuous, positive definite function $W : \Re^n \to \Re^+$ such that for every $x \in \Re^n$ the following inequality holds for all $h \in [0, \varphi(x)]$:*

$$V(x + hF(h,x)) \leq V(x) - hW(x) \tag{4.1}$$

*Then $0 \in \Re^n$ is URGAS for system (2.2). Moreover, if $W(x) := 2\sigma V(x)$, $V(x) \geq K|x|^2$ for certain $\sigma, K > 0$, then $0 \in \Re^n$ is robustly K-exponentially stable for (2.2).*

**Proof:** Notice first that by virtue of (2.3) the following claim holds:

**CLAIM:** *There exist a function $\bar{a} \in K_\infty$ such that for each $x_0 \in \Re^n$ and $h \in [0, \varphi(x_0)]$ the solution $y(t)$ of $\dot{y}(t) = F(h, x_0)$, $y(0) = x_0$ exists for all $t \in [0, h]$ and satisfies $|y(t)| \leq \bar{a}(|x_0|)$ for all $t \in [0, h]$.*

Let $R \geq 0$ (arbitrary) and consider the solution $x(t, x_0; u)$ of (2.2) corresponding to arbitrary locally bounded $u: \Re^+ \to \Re^+$ with arbitrary initial condition $x(0) = x_0$ satisfying $|x_0| \leq R$. Since $V: \Re^n \to \Re^+$ is continuous, positive definite and radially unbounded, it follows from Lemma 3.5 in [19] that there exist functions $a_1, a_2 \in K_\infty$ such that the following inequality holds:

$$a_1(|x|) \leq V(x) \leq a_2(|x|), \quad \forall x \in \Re^n \tag{4.2}$$

Using induction and (4.1) we have

$$V(x(\tau_i, x_0; u)) \leq V(x_0) \text{ for all } i \geq 0 \tag{4.3}$$

Inequality (4.3) in conjunction with (4.2) and the above claim shows that

$$|x(t, x_0; u)| \leq \bar{a}(a_1^{-1}(a_2(|x_0|))), \quad \forall t \in [0, \sup \tau_i) \tag{4.4}$$

Moreover, inequality (4.3) implies that the sequence $x(\tau_i, x_0; u)$ is bounded, which combined with the fact that $u: \Re^+ \to \Re^+$ is locally bounded, implies that $t_{\max} = \sup \tau_i = +\infty$. Consequently, estimate (4.4) guarantees Uniform Robust Lagrange Stability and Uniform Robust Lyapunov Stability. We next establish that for every $\varepsilon > 0$ it holds that

$$V(x(\tau_i, x_0; u)) \leq a_1(\bar{a}^{-1}(\varepsilon)), \text{ for all } i \in Z^+ \text{ with } \tau_i \geq \frac{a_2(R)}{w(\varepsilon, R)} \tag{4.5}$$

where

$$w(\varepsilon, R) := \min\{W(x) : a_2^{-1}(a_1(\bar{a}^{-1}(\varepsilon))) \leq |x| \leq \bar{a}(a_1^{-1}(a_2(R)))\} > 0 \tag{4.6}$$

It is clear that, by virtue of (4.2), (4.5) and the above claim it follows that $|x(t, x_0; u)| \leq \varepsilon$ for all $t \geq r + \frac{a_2(R)}{w(\varepsilon, R)}$. The previous estimate implies uniform robust global attractivity.

We next establish (4.5) by contradiction. Let $\varepsilon > 0$ (arbitrary). Suppose on the contrary that there exists $i \geq 0$ with $\tau_i \geq \frac{a_2(R)}{w(\varepsilon, R)}$ such that $V(x(\tau_i)) > a_1(\bar{a}^{-1}(\varepsilon))$. By virtue of (4.1) it follows that $V(x(\tau_k, x_0; u)) > a_1(\bar{a}^{-1}(\varepsilon))$, for all $k = 0, \ldots, i$. The previous inequality in conjunction with inequalities (4.1), (4.4) and definition (4.6) implies $V(x(\tau_{k+1}, x_0; u)) \leq V(x(\tau_k, x_0; u)) - h_k w(\varepsilon, R)$ for all $k = 0, \ldots, i-1$. Thus, we obtain $V(x(\tau_i, x_0; u)) \leq V(x_0) - w(\varepsilon, R) \sum_{k=0}^{i-1} h_k$. Notice that inequality (4.2) implies that $V(x_0) \leq a_2(R)$. Since $\tau_i = \sum_{k=0}^{i-1} h_k$, we obtain $a_1(\bar{a}^{-1}(\varepsilon)) < a_2(R) - \tau_i w(\varepsilon, R) \leq 0$, a contradiction.

Furthermore, notice that by virtue of (2.3) there exists a continuous non-decreasing mapping $\tilde{M}: \Re^+ \to \Re^+$ such that the function $\bar{a} \in K_\infty$ involved in the above claim satisfies $\bar{a}(s) \leq s\tilde{M}(s)$ (e.g., $\tilde{M}(s) = 1 + rM(s)$). If in addition we have $W(x) := 2\sigma V(x)$, $V(x) \geq K|x|^2$ for certain $\sigma, K > 0$, then (4.1) and (4.2) imply that $|x(\tau_i, x_0; u)| \leq \exp(-\sigma \tau_i) \sqrt{\frac{a_2(|x_0|)}{K}}$ for all $i \in Z^+$. By virtue of the above claim and previous inequalities we get

$$|x(t, x_0; u)| \leq \exp(-\sigma t) \exp(\sigma r) \tilde{M}\left(\sqrt{\frac{a_2(|x_0|)}{K}}\right) \sqrt{\frac{a_2(|x_0|)}{K}}$$

for all $t \in [\tau_i, \tau_{i+1})$ and $i \in Z^+$. The previous inequality implies that $0 \in \Re^n$ is robustly K-exponentially stable for (2.2). The proof is complete. ◁

The essential problem with the use of Lemma 4.1 is the knowledge of the Lyapunov function $V$. In the sequel, we will use a Lyapunov function for the continuous-time system in order to construct a Lyapunov function for its hybrid numerical approximation. To this end we use the following definition.

**Definition 4.2:** *A positive definite, radially unbounded function $V \in C^1(\Re^n; \Re^+)$ is called a Lyapunov function for system (2.1) if the following inequality holds for all $x \in (\Re^n \setminus \{0\})$:*

$$\nabla V(x) f(x) < 0 \tag{4.7}$$

In the following subsections, we show that under certain assumptions a Lyapunov function $V$ for (2.1) can be used as a Control Lyapunov Function (see [1,25,27]) in order to design the step function $\varphi : \Re^n \to (0, r]$ involved in problems (P1) and (P2). Therefore, the Lyapunov function for the original system (2.1) will be used as a Lyapunov function for its numerical approximation (2.2). The following technical results will be used in the subsequent subsections and their proofs are provided at the Appendix.

**Lemma 4.3:** Let $V \in C^1(\Re^n; \Re^+)$ be a Lyapunov function for system (2.1). Then the following statements hold:

**(i)** *There exists a locally Lipschitz, positive definite function $W : \Re^n \to \Re^+$ such that the following inequality holds for all $x \in \Re^n$:*

$$W(x) \leq -\nabla V(x) f(x) \tag{4.8}$$

**(ii)** *Let $l_f : \Re^n \to (0, +\infty)$ be a continuous function satisfying $l_f(x) \geq \sup\left\{ \frac{|f(y) - f(z)|}{|y - z|} : y, z \in \Re^n, y \neq z, \max\{V(z), V(y)\} \leq V(x) \right\}$ for all $x \in (\Re^n \setminus \{0\})$. Then for every positive constant $b > 0$, there exists a continuous, positive definite function $\widetilde{W} : \Re^n \to \Re^+$ such that the following inequality holds for all $x \in \Re^n$ and $h \in [0, \varphi(x)]$:*

$$V(z(h, x)) \leq V(x) - h \widetilde{W}(x) \tag{4.9}$$

*where*

$$\varphi(x) := \frac{b}{l_f(x)} \tag{4.10}$$

**(iii)** *Let $b > 0$, $W : \Re^n \to \Re^+$ be the function from statement (i) above and let $l_W^b : \Re^n \to \Re^+$ be a continuous positive definite function satisfying $l_W^b(x) \geq \sup\left\{ \frac{|W(y) - W(z)|}{|y - z|} : y, z \in \Re^n, y \neq z, \max\{|y|, |z|\} \leq \exp(b)|x| \right\}$ for all $x \in (\Re^n \setminus \{0\})$. If there exist constants $\varepsilon, c > 0$ such that the following inequality holds for all $x \in B_\varepsilon(0)$:*

$$|x| l_W^b(x) \leq c W(x) \tag{4.11}$$

*then for each $\lambda \in (0,1)$ inequality (4.9) holds for all $x \in \Re^n$ and $h \in [0, \varphi(x)]$ with $\widetilde{W}(x) := \lambda W(x)$ where $\varphi \in C^0(\Re^n; (0, +\infty))$ is any function satisfying*

$$\varphi(x) \leq \min\left\{ \frac{b}{l_f(x)}, \frac{(1 - \lambda) \exp(-b) W(x)}{|x| l_W^b(x) l_f(x)} \right\}, \quad \forall x \in (\Re^n \setminus \{0\}) \tag{4.12}$$

**Proposition 4.4:** *Suppose that* $f: \Re^n \to \Re^n$ *is a continuously differentiable vector field,* $0 \in \Re^n$ *is UGAS and locally exponentially stable for (2.1). Then there exist a Lyapunov function* $V \in C^1(\Re^n; \Re^+)$ *for (2.1), a symmetric, positive definite matrix* $P \in \Re^{n \times n}$ *and constants* $\varepsilon, \mu > 0$ *such that the following hold:*

$$V(x) = x'Px, \quad \forall x \in B_\varepsilon(0) \tag{4.13}$$

$$\nabla V(x) f(x) \leq -\mu |x|^2, \quad \forall x \in \Re^n \tag{4.14}$$

## 4.II. General Runge-Kutta Schemes

The following theorem provides sufficient conditions for the solvability of problem (P2) based on the Lyapunov function for the dynamical system (2.1).

**Theorem 4.5:** *Suppose that there exists an integer* $p \geq 1$ *and a Lyapunov function* $V \in C^{(p+1)}(\Re^n; \Re^+)$ *for system (2.1). Consider system (2.2) that corresponds to a Runge-Kutta scheme for (2.1) and suppose that:*

**i)** *for each fixed* $x \in \Re^n$ *the mapping* $[0, \varphi(x)] \ni h \to V(x + hF(h, x))$ *is* $(p+1)$ *times continuously differentiable.*

**ii)** *For every* $x \in \Re^n$ *and* $h \in [0, \varphi(x)]$ *there exists constant* $K > 0$ *such that* $|z(h, x) - x - hF(h, x)| \leq K h^{p+1}$, *i.e., the Runge-Kutta numerical scheme is of order* $p \geq 1$.

**iii)** *There exists a constant* $\lambda \in (0,1)$ *such that for every* $x \in \Re^n$ *it holds that* $\varphi(x) \min_{j=1,\ldots,p} K_j(x) \leq (\lambda - 1) L_f V(x)$,

*where* $\displaystyle K_j(x) := \max\left\{ \sum_{i=2}^{j} \frac{s^{i-2}}{i!} L_f^i V(x) + \frac{s^{j-1}}{(j+1)!} \frac{\partial^{j+1}}{\partial h^{j+1}} V(x + hF(h, x)) : h, s \in [0, \varphi(x)] \right\}$ *for* $j \geq 2$ *and*

$\displaystyle K_1(x) := \frac{1}{2} \max\left\{ \frac{\partial^2}{\partial h^2} V(x + hF(h, x)) : h \in [0, \varphi(x)] \right\}$.

*Then* $0 \in \Re^n$ *is URGAS for system (2.2). Moreover, if there exist constants* $\sigma, K > 0$ *such that* $\nabla V(x) f(x) \leq -2\sigma V(x)$ *and* $V(x) \geq K|x|^2$ *for all* $x \in \Re^n$ *then* $0 \in \Re^n$ *is robustly K-exponentially stable for (2.2).*

**Proof:** Since for each fixed $x \in \Re^n$ the mapping $[0, \varphi(x)] \ni h \to g(h) = V(x + hF(h, x))$ is $(p+1)$ times continuously differentiable, we have from Taylor's theorem for all $j = 1, \ldots, p$ and $h \in [0, \varphi(x)]$:

$$V(x + hF(h, x)) = g(h) \leq g(0) + \sum_{i=1}^{j} \frac{h^i}{i!} \frac{d^i g}{dh^i}(0) + \frac{h^{j+1}}{(j+1)!} \max_{0 \leq \xi \leq h} \frac{d^{j+1} g}{dh^{j+1}}(\xi) \tag{4.15}$$

Since the Runge-Kutta numerical scheme is of order $p \geq 1$, we have

$$\frac{d^i g}{dh^i}(0) = L_f^i V(x), \text{ for } i = 1, \ldots, p \tag{4.16}$$

Consequently, we obtain for all $j = 1, \ldots, p$ and $h \in [0, \varphi(x)]$:

$$V(x + hF(h, x)) \leq V(x) + h L_f V(x) + h^2 K_j(x) \tag{4.17}$$

or equivalently for all $h \in [0, \varphi(x)]$:

$$V(x+hF(h,x)) \leq V(x)+hL_f V(x)+h^2 \min_{j=1,\ldots,p} K_j(x) \qquad (4.18)$$

The inequality $\varphi(x) \min_{j=1,\ldots,p} K_j(x) \leq (\lambda-1)L_f V(x)$ in conjunction with (4.18) implies that $V(x+hF(h,x)) \leq V(x)+\lambda h L_f V(x)$. Lemma 4.1 implies that $0 \in \Re^n$ is URGAS for system (2.2). Furthermore, if there exist constants $\sigma, K > 0$ such that $\nabla V(x)f(x) \leq -2\sigma V(x)$ and $V(x) \geq K|x|^2$ for all $x \in \Re^n$ then $0 \in \Re^n$ is robustly K-exponentially stable for (2.2). The proof is complete. ◁

**Remark 4.6:**

**(a)** Theorem 4.5 shows that given a Runge-Kutta numerical scheme with order $p \geq 1$ satisfying (2.7) and a system of ODEs (2.1) with $f \in C^{(p+1)}(\Re^n; \Re^n)$ for which $0 \in \Re^n$ is GAS

"if a Lyapunov function $V \in C^{(p+1)}(\Re^n; \Re^+)$ for (2.1) is available for which there exist constants $K, \Lambda > 0$, an integer $q \geq 1$ and a neighborhood $N \subset \Re^n$ with $0 \in N$ such that $\nabla V(x)f(x) \leq -K|x|^{q+1}$ and $|f(x)| \leq \Lambda |x|^q$ for all $x \in N$, then for every $\lambda \in (0,1)$ and every compact $S \subset \Re^n$ we can find $h > 0$ sufficiently small such that $V(x+hF(h,x)) \leq V(x)+\lambda h \nabla V(x)f(x)$ for all $x \in S$"

This fact follows from (2.7) and the observation that $K_1(x) := \frac{1}{2}\max\left\{\frac{\partial^2}{\partial h^2}V(x+hF(h,x)):h \in [0,\varphi(x)]\right\} = O(|x|^{q+1})$ for $x$ close to zero. Consequently, the numerical solution of (2.1) with sufficiently small step size will give the correct dynamic behavior.

**(b)** The functions $K_j(x) := \max\left\{\sum_{i=2}^{j}\frac{s^{i-2}}{i!}L_f^i V(x)+\frac{s^{j-1}}{(j+1)!}\frac{\partial^{j+1}}{\partial h^{j+1}}V(x+hF(h,x)):h,s \in [0,\varphi(x)]\right\}$ for $j \geq 2$ and $K_1(x) := \frac{1}{2}\max\left\{\frac{\partial^2}{\partial h^2}V(x+hF(h,x)):h \in [0,\varphi(x)]\right\}$ involved in hypothesis (iii) of Theorem 4.5 are in general difficult to be computed for higher order Runge-Kutta schemes. However, for the explicit Euler scheme $F(h,x) = f(x)$ the function $K_1(x)$ can be computed without difficulty using the formula $K_1(x) := \frac{1}{2}\max\left\{f'(x)\nabla^2 V(x+hf(x))f(x):h \in [0,\varphi(x)]\right\}$. Consequently, we obtain the following corollary.

**Corollary 4.7-Explicit Euler method:** *Suppose that there exists a Lyapunov function $V \in C^2(\Re^n; \Re^+)$ for system (2.1) where $f \in C^0(\Re^n; \Re^n)$ is simply locally Lipschitz and that there exist constants $r \geq \delta > 0$, $\lambda \in (0,1)$ and a neighborhood $N \subset \Re^n$ with $0 \in N$ such that*

$$\delta q(x) \leq -2(1-\lambda)\nabla V(x)f(x), \quad \forall x \in N \qquad (4.19)$$

*where $q(x) := \max\{f'(x)\nabla^2 V(x+hf(x))f(x):h \in [0,r]\}$. Then Problem (P1) is solvable for system (2.2) with $F(h,x) := f(x)$. Particularly, $0 \in \Re^n$ is URGAS for system (2.2) if $\varphi \in C^0(\Re^n;(0,r])$ satisfies*

$$\varphi(x)q(x) \leq -2(1-\lambda)\nabla V(x)f(x), \quad \forall x \in \Re^n \qquad (4.20)$$

**Proof:** Notice that inequality (4.19) guarantees that there exists $\varphi \in C^0(\Re^n;(0,r])$ satisfying (4.20) (e.g. we may define $\varphi(x) := \delta$ for all $x \in N$, $\varphi(x) := \delta$ for all $x \notin N$ with $q(x) \leq 0$ and $\varphi(x) := \min\left\{ -\dfrac{2(1-\lambda)\nabla V(x)f(x)}{q(x)}, \delta \right\}$ for all $x \notin N$ with $q(x) > 0$). The rest is a consequence of Theorem 4.5 and the fact that $2K_1(x) \leq q(x)$ for all $x \in \mathfrak{R}^n$. ◁

**Remark 4.8:** Corollary 4.7 shows that given a system of ODEs (2.1) with $f \in C^0(\mathfrak{R}^n; \mathfrak{R}^n)$ being locally Lipschitz for which $0 \in \mathfrak{R}^n$ is GAS

> "if a Lyapunov function $V \in C^2(\mathfrak{R}^n; \mathfrak{R}^+)$ for (2.1) is available for which there exist constants $K, \Lambda > 0$, an integer $q \geq 1$ and a neighborhood $N \subset \mathfrak{R}^n$ with $0 \in N$ such that $\nabla V(x)f(x) \leq -K|x|^{2q}$ and $|f(x)| \leq \Lambda|x|^q$ for all $x \in N$, then for every $\lambda \in (0,1)$ and every compact $S \subset \mathfrak{R}^n$ we can find $h > 0$ sufficiently small such that $V(x + hf(x)) \leq V(x) + \lambda h \nabla V(x)f(x)$ for all $x \in S$"

This fact follows from (2.7) the observation that $q(x) = O(|x|^{2q})$ for $x$ close to zero. Notice the difference with Remark 4.6(a) ($\nabla V(x)f(x) \leq -K|x|^{2q}$ instead of $\nabla V(x)f(x) \leq -K|x|^{q+1}$).

The following theorem provides sufficient conditions for the solvability of problem (P2) based on the Lyapunov function for the dynamical system (2.1). It should be noticed that the hypotheses of the following theorem are of different nature from the hypotheses (i)-(iii) of Theorem 4.5. In particular, the conditions in the following theorem do not require the Lyapunov function to be smoother than $C^1$.

**Theorem 4.9:** *Consider system (2.2) that corresponds to a Runge-Kutta scheme for (2.1) of order $p \geq 1$ satisfying (2.7), (2.8a,b) for certain $\varphi \in C^0(\mathfrak{R}^n; (0,+\infty))$. Suppose that:*

**i)** *There exist a Lyapunov function $V \in C^1(\mathfrak{R}^n; \mathfrak{R}^+)$ for system (2.1) and a continuous, positive definite function $\widetilde{W} : \mathfrak{R}^n \to \mathfrak{R}^+$ such that (4.9) holds for all $x \in \mathfrak{R}^n$ and $h \in [0, \varphi(x)]$.*

**ii)** *There exists $b \geq 0$ such that $|z(h,x)| \leq \exp(b)|x|$ and $|x + hF(h,x)| \leq \exp(b)|x|$ for all $x \in \mathfrak{R}^n$ and $h \in [0, \varphi(x)]$.*

**iii)** *There exists a constant $\lambda \in (0,1)$ such that*

$$\varphi(x) \leq \left( \dfrac{(1-\lambda)\widetilde{W}(x)}{l_V^b(x)C(x)} \right)^{\dfrac{1}{p}}, \quad \forall x \in (\mathfrak{R}^n \setminus \{0\}) \tag{4.21}$$

*where $l_V^b(x) := \max\{|\nabla V(z)| : z \in \mathfrak{R}^n, |z| \leq \exp(b)|x|\}$ for all $x \in \mathfrak{R}^n$ and $C : \mathfrak{R}^n \to \mathfrak{R}^+$ is a continuous positive definite function with $|z(h,x) - x - hF(h,x)| \leq C(x)h^{p+1}$ for all $x \in \mathfrak{R}^n$ and $h \in [0, \varphi(x)]$. Then $0 \in \mathfrak{R}^n$ is URGAS for system (2.2). Moreover, if there exist constants $\sigma, K > 0$ such that $\widetilde{W}(x) \geq 2\sigma V(x)$ and $V(x) \geq K|x|^2$ for all $x \in \mathfrak{R}^n$ then $0 \in \mathfrak{R}^n$ is robustly K-exponentially stable for (2.2).*

**Proof:** Utilizing hypotheses (i) and (ii) and since $l_V^b(x) \geq \sup\left\{ \dfrac{|V(z) - V(y)|}{|z - y|} : z, y \in \mathfrak{R}^n, \max\{|y|, |z|\} \leq \exp(b)|x|, z \neq y \right\}$, we get for all $x \in (\mathfrak{R}^n \setminus \{0\})$ and $h \in [0, \varphi(x)]$:

$$V(x+hF(h,x)) \leq |V(x+hF(h,x)) - V(z(h,x))| + V(z(h,x))$$
$$\leq l_V^b(x)|x+hF(h,x) - z(h,x)| + V(x) - h\widetilde{W}(x)$$

The above inequality in conjunction with (2.8b) gives for all $x \in \Re^n$ and $h \in [0, \varphi(x)]$:

$$V(x+hF(h,x)) \leq V(x) - h\left(\widetilde{W}(x) - h^p l_V^b(x)C(x)\right)$$

where $C: \Re^n \to \Re^+$ is the continuous positive definite function with $|z(h,x) - x - hF(h,x)| \leq C(x)h^{p+1}$ for all $x \in \Re^n$ and $h \in [0, \varphi(x)]$. The above inequality in conjunction with (4.21) implies:

$$V(x+hF(h,x)) \leq V(x) - \lambda h\widetilde{W}(x), \text{ for all } x \in (\Re^n \setminus \{0\}) \text{ and } h \in [0, \varphi(x)] \tag{4.22}$$

Notice that (2.3) guarantees that (4.22) holds for $x = 0$ as well. Consequently, Lemma 4.1 implies that $0 \in \Re^n$ is URGAS for system (2.2). Moreover, if there exist constants $\sigma, K > 0$ such that $\widetilde{W}(x) \geq 2\sigma V(x)$ and $V(x) \geq K|x|^2$ for all $x \in \Re^n$ then in this case Lemma 4.1 and inequality (4.22) show that $0 \in \Re^n$ is robustly K-exponentially stable for (2.2). The proof is complete. ◁

**Remark 4.10:** The proof of Lemma 4.3 (see (A3)), inequality (2.9) and Theorem 4.8 show that given a Runge-Kutta numerical scheme with order $p \geq 1$ and a system of ODEs (2.1) with $f \in C^p(\Re^n; \Re^n)$ for which $0 \in \Re^n$ is GAS

> "if a Lyapunov function $V \in C^2(\Re^n; \Re^+)$ for (2.1) is available for which there exist constants $K, \Lambda, c > 0$, an integer $q \geq 1$ and a neighborhood $N \subset \Re^n$ with $0 \in N$ such that $\nabla V(x)f(x) \leq -K|x|^{q+1}$, $|f(x)| \leq \Lambda|x|^q$ and (4.11) with $W(x) := -\nabla V(x)f(x)$ holds for all $x \in N$, then for every $\lambda \in (0,1)$ and every compact $S \subset \Re^n$ we can find $h > 0$ sufficiently small such that $V(x+hF(h,x)) \leq V(x) + \lambda h \nabla V(x)f(x)$ for all $x \in S$"

This fact follows from (2.9) and the observation that $l_V^b(x) := \max\{|\nabla V(z)| : z \in \Re^n, |z| \leq \exp(b)|x|\} = O(|x|)$ for $x$ close to zero. The reader should notice the differences with Remark 4.6(a) (less regularity requirements). Example 4.14 below will show that in practice it is not necessary to compute $\varphi \in C^0(\Re^n; (0,+\infty))$ in order to be able to produce a qualitatively correct numerical solution; you only need to know a Lyapunov function for (2.1) with the properties described above.

The following example illustrates Remark 4.6 and Remark 4.10.

**Example 4.11:** Consider three planar systems with $\dot{x} = f_k(x)$, $k = 1,2,3$, $x = (x_1, x_2)' \in \Re^2$ with

$$f_1(x) := \begin{bmatrix} -x_1 + x_2 \\ -x_1 - x_2 \end{bmatrix}, \; f_2(x) := \begin{bmatrix} -|x|^2 x_1 + x_2 \\ -x_1 - |x|^2 x_2 \end{bmatrix}, \; f_3(x) := |x|^2 \begin{bmatrix} -x_1 + x_2 \\ -x_1 - x_2 \end{bmatrix}$$

For each of the systems we can use the Lyapunov function $V(x) = |x|^2$. Then we obtain:

$$\nabla V(x)f_1(x) = -2|x|^2, \; \nabla V(x)f_2(x) = -2|x|^4, \; \nabla V(x)f_3(x) = -2|x|^4$$

Clearly, there exist constants $\Lambda_k > 0$ ($k = 1,2,3$), integers $q_k \geq 1$ ($k = 1,2,3$) and a neighborhood $N \subset \Re^2$ with $0 \in N$ such that $|f_k(x)| \leq \Lambda_k |x|^{q_k}$ for all $x \in N$. More specifically, we have $q_1 = q_2 = 1$ and $q_3 = 3$. Remark 4.6(a) shows that for $k = 1$ and $k = 3$ we are in a position to apply any consistent Runge-Kutta numerical scheme with sufficiently small step size and produce a qualitatively correct numerical solution. The same conclusion is derived by Remark 4.10 (notice that (4.11) holds for each of the systems with $W_k(x) := -\nabla V(x) f_k(x)$, $l^b_{W_1}(x) = 4\exp(b)|x|$, $l^b_{W_2}(x) = l^b_{W_3}(x) = 8\exp(3b)|x|^3$ for a neighborhood $N \subset \Re^2$ with $0 \in N$).

On the other hand, the requirements presented in Remark 4.6(a) or Remark 4.10 are not fulfilled for $k = 2$. Similarly, the requirements presented in Remark 4.8 are not fulfilled for $k = 2$. Consequently, we cannot conclude that the application of any consistent Runge-Kutta numerical scheme with sufficiently small step size will produce a qualitatively correct numerical solution. Numerical solutions with the explicit Euler and the Heun scheme confirm these results (see Figure 1 and Figure 2). For the system $\dot{x} = f_2(x)$ both schemes applied with constant $h > 0$ exhibit an asymptotically stable limit cycle, which shrinks to the origin as $h \to 0$, but exists for all $h > 0$. This does not happen for $k = 1$ and $k = 3$: the numerical simulations with constant step size $h = 0.2$ for both numerical schemes show the correct behavior (see Figure 3 and Figure 4).

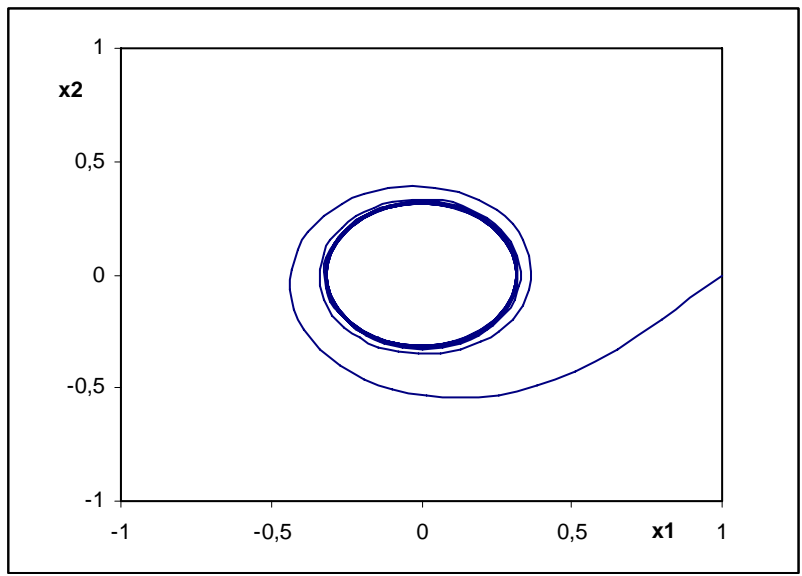

**Figure 1:** Numerical solution for $\dot{x} = f_2(x)$ with the explicit Euler method and initial condition $x = (1,0)$

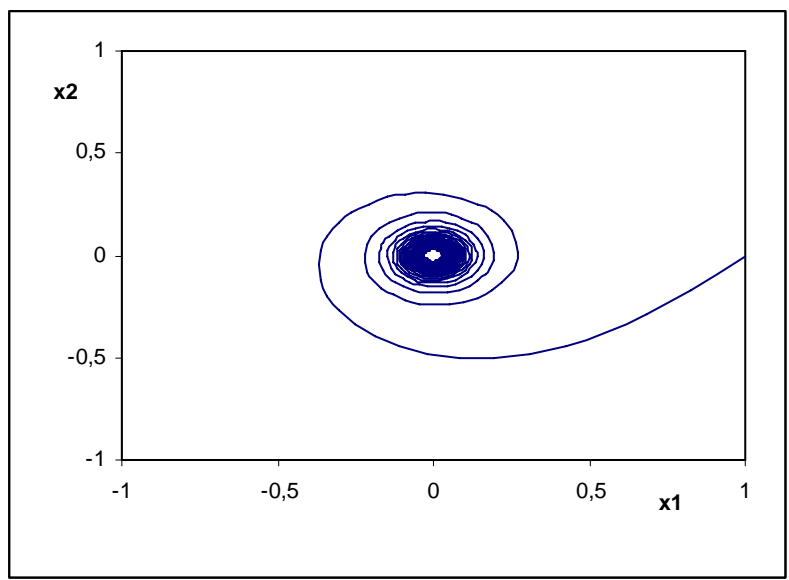

**Figure 2:** Numerical solution for $\dot{x} = f_2(x)$ with Heun's method and initial condition $x = (1,0)$

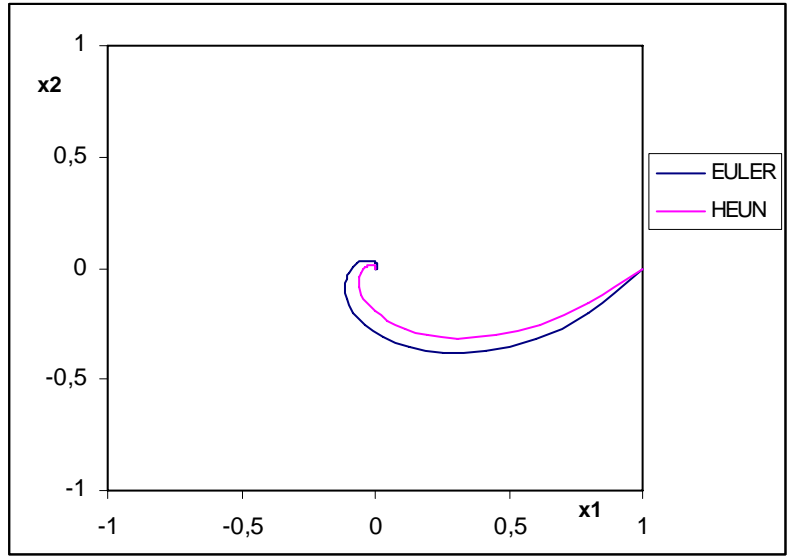

**Figure 3:** Numerical solutions for $\dot{x} = f_1(x)$ with initial condition $x = (1,0)$

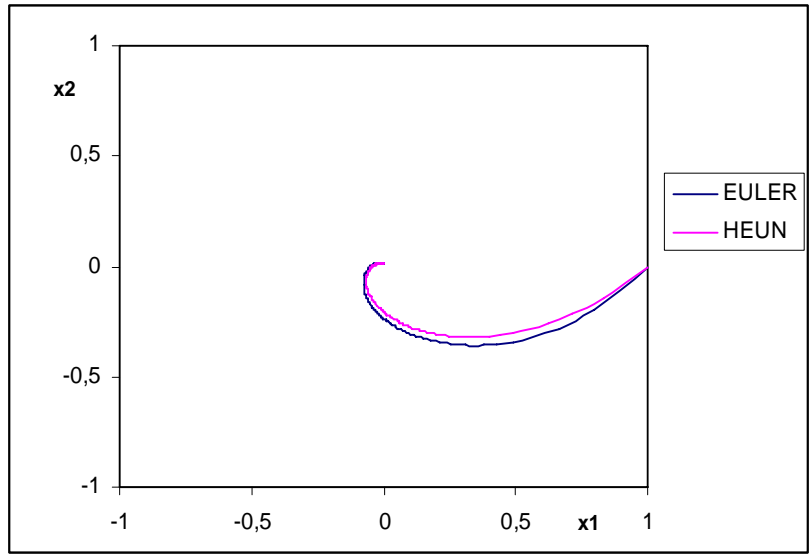

**Figure 4:** Numerical solutions for $\dot{x} = f_3(x)$ with initial condition $x = (1,0)$

We would like to emphasize that Theorem 4.5 and Theorem 4.9 do not state that there does not exist a Runge-Kutta scheme which produces an asymptotically stable approximation for system $\dot{x} = f_2(x)$, since it gives a merely sufficient but not necessary condition. In fact, for instance the implicit Euler scheme produces an asymptotically stable approximation, which we will rigorously show in Example 4.19 below. ◁

Based on the general Theorem 4.9, the following theorem shows that for the special case of a locally exponentially stable ODE system, problem (P1) is always solvable.

**Theorem 4.12:** *Consider system (2.1) and a consistent Runge-Kutta scheme with order $p \geq 1$ and $f \in C^p(\Re^n; \Re^n)$. Assume that $0 \in \Re^n$ is GAS and locally exponentially stable for (2.1). Then Problem (P1) is solvable.*

**Proof:** We are going to show that there exists $\varphi \in C^0(\Re^n; (0, +\infty))$ satisfying all requirements of Theorem 4.9.

Since $0 \in \Re^n$ is GAS and locally exponentially stable for (2.1), by virtue of Proposition 4.4, there exist a Lyapunov function $V \in C^1(\Re^n; \Re^+)$ for (2.1), a symmetric, positive definite matrix $P \in \Re^{n \times n}$ and constants $\varepsilon, \mu > 0$ such that (4.13), (4.14) hold. It follows from (4.14) that statement (i) of Lemma 4.3 holds with $W(x) := \mu |x|^2$.

Let $b > 0$. It holds that $l_W^b(x) := 2\exp(b)|x| \geq \sup\left\{ \frac{|W(y) - W(z)|}{|y - z|} : y, z \in \Re^n, y \neq z, \max\{|y|, |z|\} \leq \exp(b)|x| \right\}$ for all $x \neq 0$. Notice that (4.11) holds for all $x \in \Re^n$ with $c := 2\mu^{-1}\exp(b)$. By virtue of statement (iii) of Lemma 4.3, for each $\lambda \in (0,1)$ inequality (4.9) holds for all $x \in \Re^n$ and $h \in [0, \varphi(x)]$ with $\tilde{W}(x) := \lambda \mu |x|^2$, where $\varphi \in C^0(\Re^n; (0, +\infty))$ is any function satisfying

$$\varphi(x) \leq \frac{1}{1 + 2l_f(x)} \min\{2b, (1-\lambda)\exp(-2b)\mu\} \tag{4.23}$$

and $l_f(x) := \{|Df(z)| : z \in \Re^n, V(z) \leq V(x)\}$. Moreover, the proof of Lemma 4.3 shows that (see (A3)) for all $x \in \Re^n$ and $h \in [0, \varphi(x)]$ it holds that $|z(h, x)| \leq \exp(b)|x|$. Let $\bar{\varphi} \in C^0(\Re^n; (0, +\infty))$ the function for which (2.7), (2.8a,b) hold for all $x \in \Re^n$ and $h \in [0, \bar{\varphi}(x)]$. We notice that the inequality $|x + hF(h, x)| \leq \exp(b)|x|$ holds for all $x \in \Re^n$ and $h \in [0, \varphi(x)]$, where $\varphi \in C^0(\Re^n; (0, +\infty))$ is any function satisfying

$$\varphi(x) \leq \min\left\{ \bar{\varphi}(x), \frac{\exp(b)}{1 + \bar{\varphi}(x)M(|x|)} \right\}, \quad \forall x \in \Re^n \tag{4.24}$$

and $M : \Re^+ \to \Re^+$ is the continuous, non-decreasing function involved in (2.3).

Next we show the existence of $\varphi \in C^0(\Re^n; (0, +\infty))$ satisfying (4.21). It suffices to show that there exist constants $\delta > 0$, $\lambda \in (0,1)$ and a neighborhood $N \subset \Re^n$ with $0 \in N$ such that

$$\delta^p l_V^b(x) C(x) \leq (1 - \lambda) \mu \lambda |x|^2, \quad \forall x \in N \tag{4.25}$$

where $C : \Re^n \to \Re^+$ is a continuous positive definite function with $|z(h, x) - x - hF(h, x)| \leq C(x)h^{p+1}$ for all $x \in \Re^n$ and $h \in [0, \varphi(x)]$. Let $N = B_\rho(0)$, where $\rho := \varepsilon \exp(-b)$ and $\varepsilon > 0$ is the constant involved in (4.13). Clearly, (4.13) implies that

$$l_V^b(x) \leq 2|P|\exp(b)|x|, \quad \forall x \in N \tag{4.26}$$

where $P \in \Re^{n \times n}$ is the symmetric, positive definite matrix involved in (4.13). Notice that without loss of generality we may assume that there exists constant $K > 0$ such that (2.9) holds with $q = 1$ for all $x \in N$ and $h \in [0, \bar{\varphi}(x)]$ (the fact that there exists constant $Q > 0$ such that $|z(h, x)| \leq Q|x|$ for all $x \in N$ and $h \geq 0$ is a consequence of local exponential stability). Consequently, by virtue of (2.9), (4.26), we can guarantee that (4.25) holds for every $\lambda \in (0,1)$ with $\delta := \left( \frac{(1-\lambda)\mu\lambda}{2K|P|\exp(b)} \right)^{\frac{1}{p}}$. Therefore, from all the above we conclude that we may define

$$\varphi(x) := \min\left\{ \delta, \bar{\varphi}(x), \frac{\exp(b)}{1 + \bar{\varphi}(x)M(|x|)}, \frac{\kappa}{1 + 2l_f(x)} \right\} \quad \text{for all } x \in N \text{ and}$$

$$\varphi(x) := \min\left\{ \delta, \bar{\varphi}(x), \frac{\exp(b)}{1 + \bar{\varphi}(x)M(|x|)}, \frac{\kappa}{1 + 2l_f(x)}, \left( \frac{(1-\lambda)\lambda\mu|x|^2}{l_V^b(x)C(x)} \right)^{\frac{1}{p}} \right\} \quad \text{for all } x \notin N, \text{ where}$$

$\kappa := \min\{2b, (1-\lambda)\exp(-2b)\mu\}$, so that all requirements of Theorem 4.9 are fulfilled. The proof is complete. ◁

**Remark 4.13:** Theorem 4.12 is an existence result which does not give an explicit estimate for $\varphi(x)$, i.e., it answers (P1) but does not solve (P2). However, similar to Remark 4.8 and 4.10 we can conclude that the numerical approximation is asymptotically stable on each compact set S for sufficiently small step size h. Note that local exponential stability is not a necessary condition for asymptotic stability of explicit Runge-Kutta schemes, as Example 4.11 shows ($0 \in \Re^2$ is GAS but not locally exponentially stable for system $\dot{x} = f_3(x)$).

Due to the complexity of calculations needed in order to satisfy the requirements of Theorem 4.5 and Theorem 4.9, we suggest the following algorithm for Runge-Kutta schemes and a given $\lambda \in (0,1)$:

"If $V(x + hF(h,x)) \leq V(x) + \lambda h L_f V(x)$ then the time step $h > 0$ is allowed.

Otherwise halve the step size (i.e., $h = h/2$) and check again"

where $V \in C^2(\Re^n; \Re^+)$ is a Lyapunov function for (2.1) for which there exist a constant $K > 0$ and a neighborhood $N \subset \Re^n$ with $0 \in N$ such that $\nabla V(x) f(x) \leq -K|x|^2$ for all $x \in N$. Therefore, in practice we do not have to compute the step function $\varphi(x)$ that guarantees robust global asymptotic stability of the numerical approximation. The following example illustrates this point.

**Example 4.14:** In this example we consider four different explicit numerical schemes: the explicit Euler scheme, Heun's scheme, the improved polygonal scheme (a second order method) and Kutta's 3$^{rd}$ order scheme. The numerical schemes are applied to the planar system:

$$\dot{x}_1 = -x_1 + x_2^2 \qquad \dot{x}_2 = -x_2 - x_1 x_2 \qquad (4.27)$$

using the Lyapunov function $V(x) = (x_1^2 + x_2^2)/2$. For all numerical schemes (except the explicit Euler method) the calculation of the maximum allowed time-step by using Theorem 4.5 or Theorem 4.9 is very complicated. However, we have determined the maximum $h > 0$ for which the inequality $V(x + hF(h,x)) \leq V(x) + \lambda h L_f V(x)$ with $\lambda = \frac{1}{2}$ holds. The following figure shows the graph of the maximum allowable time-step for the four numerical methods with $x = (x_1, 1)' \in \Re^2$ and varying $x_1 \in \Re$.

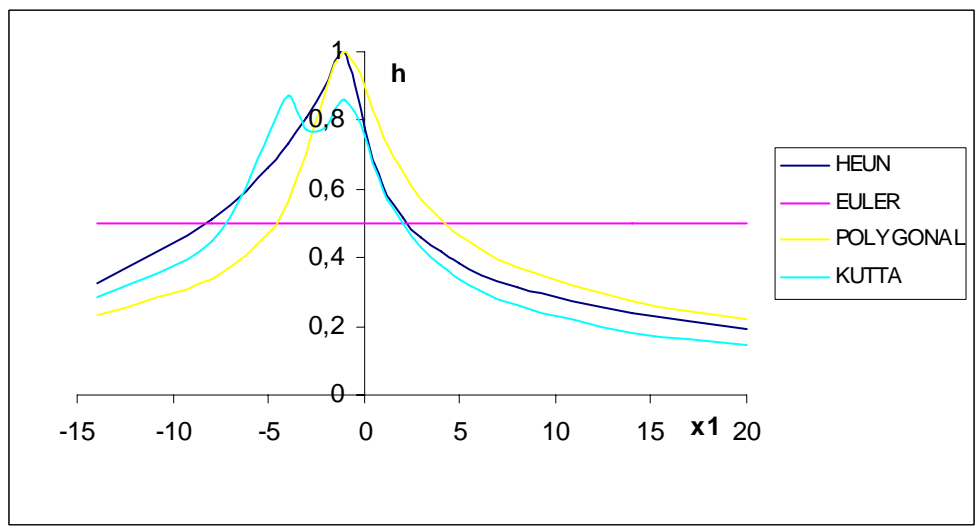

**Figure 5:** Maximum allowable time step for the explicit Euler method, Heun's method, the improved polygonal method and Kutta's 3$^{rd}$ order method for (4.27) with $x = (x_1, 1)' \in \Re^2$

It should be noticed that all higher order schemes allow greater time-steps for $x_1$ close to zero than the time-steps allowed by the explicit Euler method (notice that for $x = (x_1, 1)' \in \Re^2$ and $\lambda = \frac{1}{2}$ the maximum allowable time-step for which the inequality $V(x+hF(h,x)) \leq V(x) + \lambda h L_f V(x)$ holds for the explicit Euler method is $h = \frac{1}{2}$). However, for large values of $|x_1|$ the maximum allowable time-step for higher order schemes are considerably smaller than the time-step allowed by the explicit Euler method. It is clear that a higher order method does not necessarily allow higher values for the maximum allowable time step for which the inequality $V(x+hF(h,x)) \leq V(x) + \lambda h L_f V(x)$ holds. ◁

## 4.III. Implicit Runge-Kutta Schemes

In this subsection we show how Lyapunov function based arguments can be used for implicit schemes. In order to keep the presentation technically simple, we restrict ourselves to the implicit Euler scheme for which we can prove the following result.

**Theorem 4.15-Implicit Euler Method:** *Suppose that there exists a convex Lyapunov function for (2.1), where $f \in C^0(\Re^n; \Re^n)$ is simply locally Lipschitz. Let $\overline{\varphi} \in C^0(\Re^n; (0,+\infty))$ be such that the equation $Y = x + hf(Y)$ has a unique solution $Y \in \Re^n$ for all $h \in [0, \overline{\varphi}(x)]$ and $x \in \Re^n$. Then for each $r > 0$, $0 \in \Re^n$ is URGAS for the corresponding system (2.2) with $F(h,x) := f(Y)$, $\varphi(x) := \min\{\overline{\varphi}(x), r\}$, where $Y = x + hf(Y)$.*

**Proof:** Define the functions

$$W_1(x) := \min\left\{-\nabla V(y) f(y) : y \in \Re^n, V(x) \geq V(y) \geq \frac{1}{2} V(x)\right\}, \quad W_2(x) := \frac{1}{2r} V(x) \tag{4.28}$$

By virtue of (4.7) both functions are continuous and positive definite. Since $V \in C^1(\Re^n; \Re^+)$ is convex the following inequality holds for all $x_1, x_2 \in \Re^n$:

$$V(x_1) + \nabla V(x_1) x_2 \leq V(x_1 + x_2) \tag{4.29}$$

We apply (4.29) with $x_1 = Y$ and $x_2 = -hf(Y)$, where $Y = x + hf(Y)$ and $h \in [0, \overline{\varphi}(x)]$. We get:

$$V(x) = V(Y - hf(Y)) \geq V(Y) - h \nabla V(Y) f(Y) \tag{4.30}$$

By virtue of (4.7), (4.30) implies that $V(Y) \leq V(x)$. We distinguish the following cases:

<u>Case 1:</u> $V(Y) \geq \frac{1}{2} V(x)$. In this case (4.30) in conjunction with definition (4.28) of $W_1$ we obtain:

$$V(Y) + hW_1(x) \leq V(x) \tag{4.31}$$

<u>Case 2:</u> $V(Y) < \frac{1}{2} V(x)$. In this case definition (4.28) of $W_2$ implies for all $h \in [0, r]$:

$$V(Y) + hW_2(x) \leq V(x) \tag{4.32}$$

Consequently, in any case we obtain for all $h \in [0, \varphi(x)]$ and $x \in \Re^n$:

$$V(Y) \leq V(x) - hW(x) \tag{4.33}$$

where $Y = x + hf(Y)$ and $W(x) := \min\{W_1(x), W_2(x)\}$ is a positive definite function. Lemma 4.1 yields the assertion. The proof is complete. ◁

The following corollary shows that Theorem 4.15 can be seen as a nonlinear generalization of the well-known A-stability property of the implicit Euler method.

**Corollary 4.16:** Consider the system of ODEs $\dot{x} = Ax$, $x \in \Re^n$ where $A \in \Re^{n \times n}$ is a Hurwitz matrix. Then the implicit Euler method is URGAS for arbitrary step size $h > 0$.

**Proof:** As pointed out in Section 2, the implicit Euler method is well defined for each step size $h > 0$. Furthermore, the system $\dot{x} = Ax$, $x \in \Re^n$ admits the quadratic Lyapunov function $V(x) = x'Px$, where $P \in \Re^{n \times n}$ is a symmetric, positive definite matrix (see Theorem 5.7.18 in [27]). This Lyapunov function is obviously convex and Theorem 4.15 yields the assertion. ◁

**Remark 4.17:** The main result in [9] guarantees that if $n \neq 4,5$ then there exists a homeomorphism $\Phi : \Re^n \to \Re^n$ with $\Phi(0) = 0$, being a diffeomorphism on $(\Re^n \setminus \{0\})$ and $C^1$ on $\Re^n$ such that the transformed system (2.1) $\dot{y} = D\Phi(\Phi^{-1}(y))f(\Phi^{-1}(y))$ admits the convex Lyapunov function $V(y) := \frac{1}{2}|y|^2$. Consequently, the implicit Euler can be applied to the transformed system. However, for numerical purposes, the method is not practical, since the homeomorphism $\Phi : \Re^n \to \Re^n$ is usually not available. On the other hand, for certain classes of systems Theorem 4.15 is directly applicable. One such class are the so called gradient systems, as shown in the following example.

**Example 4.18:** Consider the following class of systems:

$$\dot{x} = f(x) = -(P(x) + G(x))(\nabla V(x))', \quad x \in \Re^n \tag{4.34}$$

where $V \in C^2(\Re^n; \Re^+)$ is a positive definite, radially unbounded function with positive definite Hessian and $\nabla V(0) = 0$, $P(x) \in \Re^{n \times n}$ is a symmetric positive definite matrix with locally Lipschitz elements and $G(x) \in \Re^{n \times n}$ is a matrix with locally Lipschitz elements with $G'(x) = -G(x)$ for all $x \in \Re^n$. The class of systems of the form (4.34) is a generalization of the class of the so-called gradient systems (see [28]).

It follows from Theorem 4.15 that the implicit Euler scheme can be applied successfully for the numerical approximation of the solutions of (4.34) for every $r > 0$, $\lambda \in (0,1)$ with $\varphi(x) := \min\left\{\frac{\lambda}{L_\lambda(x) + \gamma(x)}, r\right\}$, where $\gamma : \Re^n \to \Re^+$ is a continuous function with $|f(x)| \leq |x|\gamma(x)$ for all $x \in \Re^n$, $L_\lambda : \Re^n \to (0, +\infty)$ is a continuous function with $L_\lambda(x) \geq \sup\left\{\frac{|f(z) - f(y)|}{|z - y|} : z, y \in B_\lambda(x), z \neq y\right\}$ for all $x \in (\Re^n \setminus \{0\})$ and $B_\lambda(x) := \{y \in \Re^n : |y - x| \leq \lambda|x|\}$. ◁

**Example 4.19:** Consider again systems $\dot{x} = f_k(x)$, $k = 1,2,3$ of Example 4.11. Since these systems admit the convex Lyapunov function $V(x) = |x|^2$, it follows that the implicit Euler scheme produces asymptotically stable solutions for all systems $\dot{x} = f_k(x)$, $k = 1,2,3$ of Example 4.11. ◁

## 5. Application of the Small-Gain Step Selection Methodology

Consider the infinite-dimensional dynamical system:

$$\frac{\partial x}{\partial t}(t,z) + c\frac{\partial x}{\partial z}(t,z) = b(x(t,z))x(t,z), \ z \in (0,1] \quad (5.1)$$
$$x(t,0) = 0$$

with $x(t,z) \in \Re$, $b: \Re \to \Re$ being locally Lipschitz, $c > 0$ and initial condition $x(0,z) = x_0(z)$, where $x_0 \in C^1([0,1];\Re)$ with $x_0(0) = \frac{dx_0}{dz}(0) = 0$, under the following hypothesis

**(H)** There exists constant $K \geq 0$ such that $b(x) \leq K$ for all $x \in \Re$.

Using the method of characteristics and hypothesis (H), it can be shown that the infinite-dimensional dynamical system (5.1) admits a unique classical solution $x(t,\cdot) \in C^1([0,1];\Re)$ with $x(t,0) = \frac{\partial x}{\partial z}(t,0) = 0$ for all $t \geq 0$. Moreover, the zero solution is globally asymptotically stable, since for every $x_0 \in C^1([0,1];\Re)$ with $x_0(0) = \frac{dx_0}{dz}(0) = 0$, the unique classical solution $x(t,\cdot) \in C^1([0,1];\Re)$ of (5.1) with initial condition $x(0,z) = x_0(z)$ satisfies $x(t,z) = 0$ for all $t \geq c^{-1}z$ (uniform global attractivity).

Using a uniform space grid of $n+1$ points with space discretization step $\Delta z = \frac{1}{n}$, setting $x_i(t) = x(t,i\Delta z)$, $i = 0,1,...,n$ and approximating the spatial derivative by the backward difference scheme $\frac{\partial x}{\partial z}(t,i\Delta z) \approx \frac{x(t,i\Delta z) - x(t,(i-1)\Delta z)}{\Delta z} = \frac{x_i(t) - x_{i-1}(t)}{\Delta z}$ for $i = 1,...,n$, we obtain the following set of ordinary differential equations:

$$\dot{x}_1 = -\left(\frac{c}{\Delta z} - b(x_1)\right)x_1$$
$$\dot{x}_i = -\left(\frac{c}{\Delta z} - b(x_i)\right)x_i + \frac{c}{\Delta z}x_{i-1}, \ i = 2,...,n \quad (5.2)$$

It is clear that system (5.2) has the structure of system (3.1) with $a_i(x_i) = \frac{c}{\Delta z} - b(x_i)$ for $i = 1,...,n$. Moreover, if the space discretization step is selected so that

$$K\Delta z < c \quad (5.3)$$

where $K \geq 0$ is the constant involved in Hypothesis (H), then inequalities (3.2) hold as well with $L_i = \frac{c}{\Delta z} - K$ for $i = 1,...,n$. Theorem 3.1 allows us to conclude that for every $h > 0$ the numerical scheme:

$$x_1(t+h) = \frac{x_1(t)}{1 + h\left(\frac{c}{\Delta z} - b(x_1(t))\right)}$$

$$x_i(t+h) = \frac{x_i(t) + \frac{ch}{\Delta z}x_{i-1}(t)}{1 + h\left(\frac{c}{\Delta z} - b(x_i(t))\right)}, \ i = 2,...,n \quad (5.4)$$

will give the correct qualitative behavior. The reader should notice that for the case $b(x) \equiv 0$ inequality (5.3) is automatically satisfied (with $K = 0$) and the numerical scheme (5.4) is related to the so-called implicit upwind numerical scheme for the advection equation, which is unconditionally stable.

# 6. Applications of the Lyapunov-based Step Selection Methodology

In this section we present some applications for the Lyapunov-based step selection methodology that was provided in the Section 4. It should be emphasized that the Lyapunov-based step selection methodology can be (in principle) applied to all dynamical systems for which a Lyapunov function is known with a globally asymptotically stable and locally exponentially stable equilibrium (Theorem 4.7). However, as the following applications show there are certain classes of systems that we can guarantee more properties or require special attention.

Application 1: Solution of Nonlinear Programming Problems

There are many nonlinear programming problems which can be solved by constructing a dynamical system with a globally asymptotically stable equilibrium point which coincides with the minimizer of the nonlinear programming problem (see [5,29,30,31]). A special feature for such methods is that a Lyapunov function is available; however the position of the equilibrium point is not known (this is what we seek). Consider the following nonlinear programming problem:

$$\min f(x), x \in \Re^n$$
$$\text{s.t.} \quad \text{(P)}$$
$$Ax = b$$

where $f \in C^3(\Re^n; \Re)$ is strictly convex and radially unbounded with positive definite Hessian and $A \in \Re^{m \times n}$, $b \in \Re^m$ with $m < n$ satisfies $\det(AA') \neq 0$. Under the previous hypotheses there is a global minimum $x^* \in \Re^n$ of problem (P). Moreover, there exists a vector $z^* \in \Re^m$ such that $(x^*, z^*) \in \Re^{n+m}$ is the unique solution of the equations:

$$\nabla f(x) + z'A = 0$$
$$Ax = b \tag{6.1}$$

Problem (P) may be solved by means of differential equations if we further assume that the function $G(x) = \left| \nabla f(x)(I - A'(AA')^{-1}A) \right|^2 + |Ax - b|^2$ is radially unbounded. Indeed, the system

$$\dot{x} = -\left( \nabla^2 f(x) (\nabla f(x) + z'A)' + A'(Ax - b) \right)$$
$$\dot{z} = -A(\nabla f(x) + z'A)' \tag{6.2}$$

has the unique equilibrium point $(x^*, z^*) \in \Re^{n+m}$, which is UGAS for (6.2). This fact can be proved by using the Lyapunov function $V(x, z) = \frac{1}{2} |\nabla f(x) + z'A|^2 + \frac{1}{2} |Ax - b|^2$ (notice that it is radially unbounded). Notice that $\dot{V} = -\left( |\dot{x}|^2 + |\dot{z}|^2 \right)$ for all $(x, z) \in \Re^{n+m}$. Thus the dynamical system (6.2) can be solved by means of Runge-Kutta methods with a Lyapunov-based step selection methodology: each Runge-Kutta method applied to the dynamical system (6.2) will yield a method for the solution of the nonlinear programming problem (P).

Here we will discuss the explicit Euler method. Indeed, the requirements of Corollary 4.7 are fulfilled. To see this let $r > 0$, $\lambda \in (0,1)$ and notice that the function $q: \Re^{n+m} \to (0, +\infty)$ involved in (4.19), (4.20) satisfies

$$q(x, z) \leq |(\dot{x}, \dot{z})|^2 p(x, z) \tag{6.3}$$

where $p(x, z) := \max \left\{ \left| \nabla^2 V(y, \xi) \right| : |(y - x, \xi - z)| \leq r |(\dot{x}, \dot{z})| \right\}$ is a continuous function which can be evaluated without knowledge of the equilibrium point $(x^*, z^*) \in \Re^{n+m}$. Let $N \subset \Re^{n+m}$ be defined by $N := \{(x, z) \in \Re^{n+m} : |(x - x^*, z - z^*)| < c\}$, where $c > 0$ is any positive constant. Then condition (4.19) is implied by the following inequality:

$$\delta\, p(x,z) \leq 2(1-\lambda), \quad \forall (x,z) \in N \tag{6.4}$$

and it is clear that (6.4) holds with $\delta > 0$ sufficiently small. Notice that inequality (4.20) is satisfied with $\varphi(x,z) \leq \min\left\{\dfrac{2(1-\lambda)}{p(x,z)}; r\right\}$. Consequently, Corollary 4.7 guarantees that for every $(x_0, z_0) \in \Re^{n+m}$, the sequence $\{(x_k, z_k) \in \Re^{n+m}\}_0^\infty$ generated by the recursive formulae

$$\begin{aligned}
x_{k+1} &= x_k - h_k\left(\nabla^2 f(x_k)(\nabla f(x_k) + z_k' A)' + A'(Ax_k - b)\right)\\
z_{k+1} &= z_k - h_k A(\nabla f(x_k) + z_k' A)'
\end{aligned} \tag{6.5}$$

will achieve convergence to the (unknown) equilibrium point $(x^*, z^*) \in \Re^{n+m}$ of (6.2), provided that (the discretization step size) $h_k > 0$ satisfies $h_k \leq \min\left\{\dfrac{2(1-\lambda)}{p(x_k, z_k)}; r\right\}$.

Application 2: Control Systems under Feedback Control

One class of dynamical systems for which a Lyapunov function is known is the class of control systems for which a continuous feedback stabilizer is designed by using a Lyapunov-based methodology (see [1,3,19,25]). This is evident for the class of so-called "triangular control systems" (see [3]). Consider the following triangular control system:

$$\begin{aligned}
\dot{x}_i &= f_i(x_1, \ldots, x_i) + g_i(x_1, \ldots, x_i)x_{i+1}, \quad i=1,\ldots,n-1\\
\dot{x}_n &= f_n(x) + g_n(x)u\\
x &= (x_1, \ldots, x_n)' \in \Re^n, u \in \Re
\end{aligned} \tag{6.6}$$

where $f_i : \Re^i \to \Re$, $g_i : \Re^i \to \Re$ ($i=1,\ldots,n$) are locally Lipschitz functions with $f_i(0) = 0$ ($i=1,\ldots,n$) and $g_i(y) > 0$ for all $y \in \Re^i$ ($i=1,\ldots,n$).

Using backstepping (see [3]), we are in a position to construct a smooth function $k : \Re^n \to \Re$ with $k(0) = 0$ and a positive definite and radially unbounded smooth function $V : \Re^n \to \Re^+$ such that

$$\nabla V(x) \begin{bmatrix} f_1(x_1) + g_1(x_1)x_2 \\ \vdots \\ f_n(x) + g_n(x)k(x) \end{bmatrix} \leq -\sigma V(x), \quad \forall x \in \Re^n \tag{6.7}$$

for certain constant $\sigma > 0$. Moreover, $0 \in \Re^n$ is locally exponentially stable for the closed-loop system (6.6) with $u = k(x)$ and for every $\Delta \geq 0$ there exist constants $K_1, K_2 > 0$ such that the following inequality holds:

$$K_1 |x|^2 \leq V(x) \leq K_2 |x|^2, \text{ for all } x \in \Re^n \text{ with } |x| \leq \Delta \tag{6.8}$$

Consequently, Corollary 4.7 guarantees that the explicit Euler method can be used for the numerical approximation of the closed-loop system (6.6) with $u = k(x)$. Furthermore, Corollary 4.7 can be used in order to obtain explicit estimate of the allowable discretization time step for the explicit Euler method. Indeed, notice that all requirements of Corollary 4.7 hold (notice that (6.7), in conjunction with (6.8) show that $\nabla V(x) \begin{bmatrix} f_1(x_1) + g_1(x_1)x_2 \\ \vdots \\ f_n(x) + g_n(x)k(x) \end{bmatrix} \leq -c|x|^2$ for appropriate $c > 0$ for every bounded neighborhood of the origin). Formula (4.20) (combined with (6.7)) provides an explicit upper bound for the function $\varphi \in C^0(\Re^n; (0, r])$:

$$\varphi(x) \leq \min\left\{-\frac{2(1-\lambda)\sigma V(x)}{p(x)|F(x)|^2}; r\right\}, \quad \forall x \in (\Re^n \setminus \{0\}) \tag{6.9}$$

where $F(x) := \begin{bmatrix} f_1(x_1) + g_1(x_1)x_2 \\ \vdots \\ f_n(x) + g_n(x)k(x) \end{bmatrix}$ and $p: \Re^n \to (0, +\infty)$ is defined by:

$$p(x) := \max\left\{\left|\nabla^2 V(y)\right| : |y - x| \leq r|F(x)|\right\} \tag{6.10}$$

Other Runge-Kutta numerical schemes can be used as well. Notice that the backstepping procedure achieves the construction of the Lyapunov function $V: \Re^n \to \Re^+$ by constructing a diffeomorphism $\Phi: \Re^n \to \Re^n$ is with $\Phi(0) = 0$ such that $V(x) = \Phi'(x)P\Phi(x)$, where $P \in \Re^{n \times n}$ is a symmetric, positive definite matrix. Then Theorem 4.15 guarantees that the implicit Euler can be used as well for the transformed closed-loop system (6.6) with $u = k(x)$:

$$\dot{z} = \widetilde{F}(z) := D\Phi(x)F(x)\big|_{x=\Phi^{-1}(z)}, \text{ with } F(x) := \begin{bmatrix} f_1(x_1) + g_1(x_1)x_2 \\ \vdots \\ f_n(x) + g_n(x)k(x) \end{bmatrix} \tag{6.11}$$

It follows that for every $r > 0$, $\lambda \in (0,1)$, the implicit Euler scheme can be applied for (6.11) with $\varphi(z) := \min\left\{\frac{\lambda}{L_\lambda(z) + \gamma(z)}, r\right\}$, where $\gamma: \Re^n \to \Re^+$ is a continuous function with $|\widetilde{F}(z)| \leq |z|\gamma(z)$ for all $z \in \Re^n$, $L_\lambda: \Re^n \to (0, +\infty)$ is a continuous function with $L_\lambda(z) \geq \sup\left\{\frac{|\widetilde{F}(x) - \widetilde{F}(y)|}{|x - y|} : x, y \in B_\lambda(z), x \neq y\right\}$ for all $z \in (\Re^n \setminus \{0\})$ and $B_\lambda(z) := \left\{y \in \Re^n : |y - z| \leq \lambda|z|\right\}$. This fact was observed in [17].

Application 3: Explicit Methods for Stiff Linear Systems

Even for linear stiff systems the results provided by Theorems 4.5, 4.9 and 4.12 have important consequences. Consider the following linear system:

$$\dot{x} = Ax \quad ; \quad x = (x_1, \ldots, x_n)' \in \Re^n \tag{6.12}$$

where $A \in \Re^{n \times n}$ is a diagonalizable Hurwitz matrix with eigenvalues $\lambda_1, \ldots, \lambda_n \in C$ (all of them with negative real part). The standard criterion used in numerical analysis for the stability of a Runge-Kutta scheme requires that for all $i = 1, \ldots, n$, the complex number $h\lambda_i$ lies inside the region $S = \{z \in C : |R(z)| \leq 1\}$, where $R(z)$ is the stability function of the scheme and $h$ is the (constant) discretization step size. The possibility of using larger discretization step size for explicit Runge-Kutta methods than the one allowed by the classical analysis must be considered. This possibility was considered in [2] where it was shown that a sequence of "small" time steps can allow a "big" time step for explicit ODE solvers.

Here for simplicity, we consider the Explicit Euler scheme. The fact that $A \in \Re^{n \times n}$ is a Hurwitz matrix guarantees the existence of a symmetric positive definite matrix $P \in \Re^{n \times n}$ so that $PA + A'P < 0$. Then Corollary 4.7 implies that for every $\lambda \in (0,1)$, $r > 0$ the step-size function $\varphi \in C^0(\Re^n; (0, r])$ satisfying the inequality

$$\varphi(x) \leq \min\left\{-\frac{(1-\lambda)x'(A'P + PA)x}{x'A'PAx}; r\right\}, \quad \forall x \in (\Re^n \setminus \{0\}) \tag{6.13}$$

guarantees that the numerical solution produced by the explicit Euler scheme has the correct qualitative behavior. Notice that the quantity $-\dfrac{x'(A'P+PA)x}{x'A'PAx}$ depends heavily on the direction of the vector $x \in \Re^n$ and can allow greater discretization step sizes than the one produced by classical stability analysis.

**Example 5.1:** Consider the singularly perturbed and stiff linear system

$$\dot{x}_1 = -1000 x_1 \quad \dot{x}_2 = x_1 - x_2 \qquad (6.14)$$

Classical analysis (Dahlquist test) for the explicit Euler scheme with constant step size requires that the discretization step must satisfy $h \leq 1/500$.

Let $P := \dfrac{1}{2} I$. We get $\varphi(x) \leq \min\left\{ 2(1-\lambda) \dfrac{1000 x_1^2 + x_2^2 - x_1 x_2}{1000^2 x_1^2 + (x_1 - x_2)^2} \,;\, r \right\}$. Notice that for $x_1 = 0$ we have $\varphi(x) \leq \min\{ 2(1-\lambda) \,;\, r \}$ (a very large value for the discretization step size) and for $x_2 = 0$ we have $\varphi(x) \leq \min\left\{ 2(1-\lambda) \dfrac{1000}{1000^2 + 1} \,;\, r \right\}$ (a value for the discretization step size comparable to the value provided by the limit of the classical analysis).

Figure 6 shows the numerical solution obtained by using the explicit Euler scheme for (6.14) with step size $h = \varphi(x) = \min\left\{ 2(1-\lambda) \dfrac{1000 x_1^2 + x_2^2 - x_1 x_2}{1000^2 x_1^2 + (x_1 - x_2)^2} \,;\, r \right\}$, $\lambda = 0.6$, $r = 1$. Figure 7 shows the sequence of discretization time steps for the explicit Euler scheme with step size $h = \varphi(x) = \min\left\{ 2(1-\lambda) \dfrac{1000 x_1^2 + x_2^2 - x_1 x_2}{1000^2 x_1^2 + (x_1 - x_2)^2} \,;\, r \right\}$, $\lambda = 0.6$, $r = 1$. The initial conditions were $x_1(0) = 1$, $x_2(0) = 1.1$. It can be seen that the numerical solution presents the correct qualitative behavior. Moreover, Figure 7 shows a repeated pattern for the step selection: after a number of "small" time-steps the scheme allows a "big" time step. This is in agreement with the results in [2]. However, the "big" time step creates errors in the evaluation of the state variable $x_1(t)$ and this explains the behavior presented in Figure 6. After 500 explicit Euler steps the value of time is $t = 12.71372$: if we had applied 500 Euler steps with constant step size we would have reached at most $t = 1$.

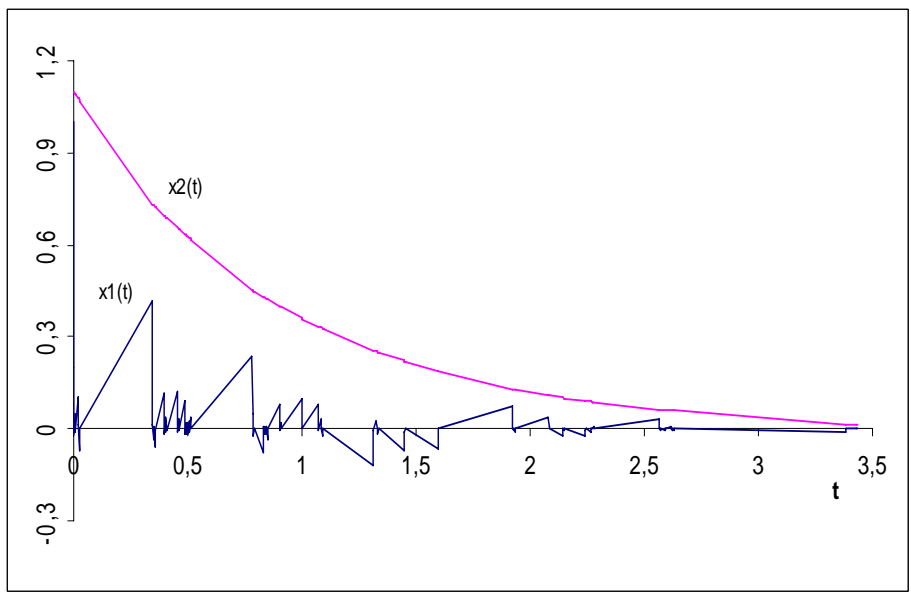

**Figure 6:** Numerical solution produced by the Explicit Euler scheme for (5.14) with

$$h = \varphi(x) = \min\left\{ 2(1-\lambda) \dfrac{1000 x_1^2 + x_2^2 - x_1 x_2}{1000^2 x_1^2 + (x_1 - x_2)^2} \,;\, r \right\}, \lambda = 0.6 , r = 1$$

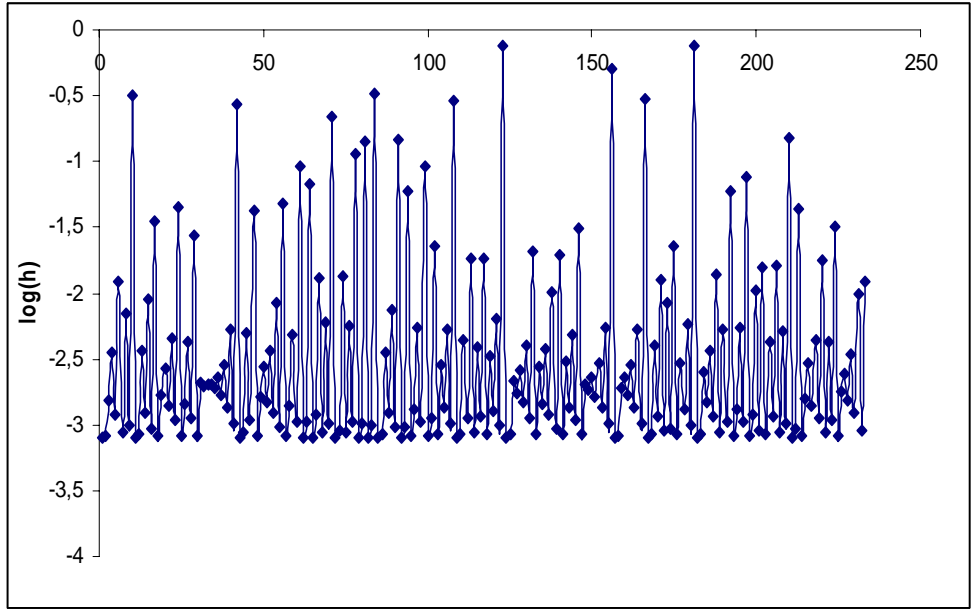

**Figure 7:** The sequence of time steps for the Explicit Euler scheme for (5.14) with
$$h = \varphi(x) = \min\left\{ 2(1-\lambda)\frac{1000x_1^2 + x_2^2 - x_1 x_2}{1000^2 x_1^2 + (x_1 - x_2)^2} ; r \right\}, \; \lambda = 0.6, \; r = 1$$

The accuracy of the numerical solution can be improved by imposing a higher value for $\lambda$. Figure 8 shows the numerical solution obtained by using the explicit Euler scheme for (5.14) with step size $h = \varphi(x) = \min\left\{ 2(1-\lambda)\frac{1000x_1^2 + x_2^2 - x_1 x_2}{1000^2 x_1^2 + (x_1 - x_2)^2} ; r \right\}$, $\lambda = 0.9$, $r = 1$. Figure 9 shows the sequence of discretization time steps for the explicit Euler scheme with step size $h = \varphi(x) = \min\left\{ 2(1-\lambda)\frac{1000x_1^2 + x_2^2 - x_1 x_2}{1000^2 x_1^2 + (x_1 - x_2)^2} ; r \right\}$, $\lambda = 0.9$, $r = 1$. Indeed, for $\lambda = 0.9$ "smaller" time steps are imposed and the accuracy is improved (compare Figure 8 with Figure 6). However, for $\lambda = 0.9$ after 500 explicit Euler steps the value of time is $t = 3.798454$.

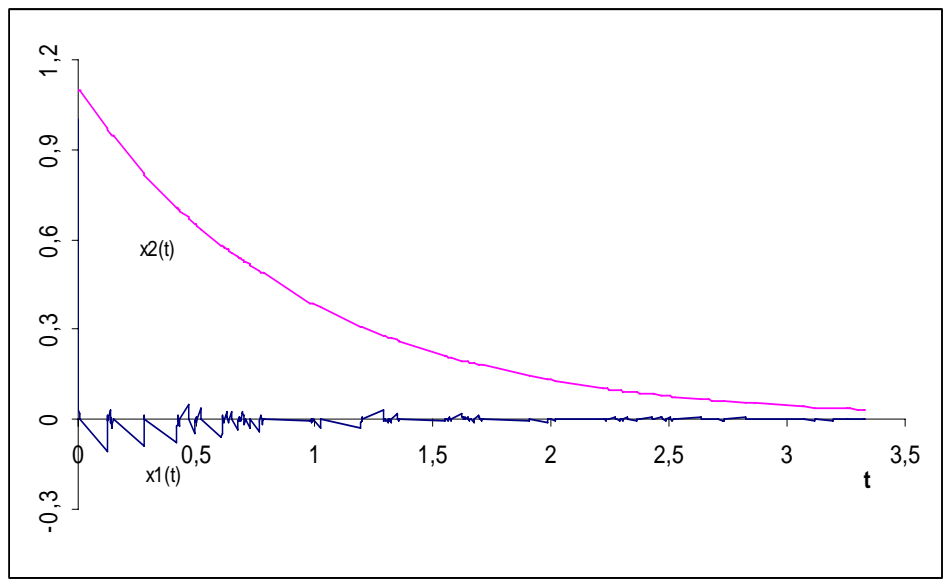

**Figure 8:** Numerical solution produced by the Explicit Euler scheme for (5.14) with
$$h = \varphi(x) = \min\left\{ 2(1-\lambda)\frac{1000x_1^2 + x_2^2 - x_1 x_2}{1000^2 x_1^2 + (x_1 - x_2)^2} ; r \right\}, \; \lambda = 0.9, \; r = 1$$

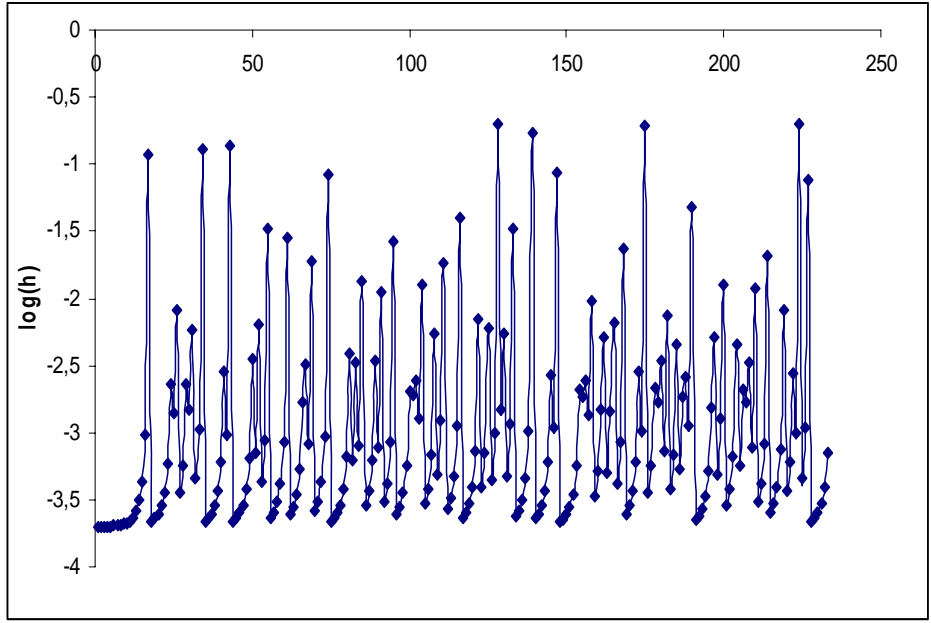

**Figure 9:** The sequence of time steps for the Explicit Euler scheme for (5.14) with

$$h = \varphi(x) = \min\left\{ 2(1-\lambda)\frac{1000 x_1^2 + x_2^2 - x_1 x_2}{1000^2 x_1^2 + (x_1 - x_2)^2} ; r \right\}, \; \lambda = 0.9, \; r = 1$$

It is clear that a trade-off between the allowable time steps and the accuracy of the numerical solution exists, as expected. ◁

## 7. Stabilization of Numerical Schemes and Control of the Global Discretization Error

It is well-known that the numerical solution produced by a Runge-Kutta numerical scheme of order $p \geq 1$ results to global discretization error $O(h^p)$ for a finite time interval, where $h > 0$ is the maximum applied discretization time step. In this section we prove a similar result for the infinite interval $[0, +\infty)$, on which, however, we will get a lower order of convergence than $p \geq 1$.

Suppose that $0 \in \Re^n$ is K-exponentially stable for (2.1), i.e., there exist a constant $\sigma > 0$ and a function $a \in K_\infty$ such that for every $z_0 \in \Re^n$ the solution $z(t, z_0)$ of (2.1) with initial condition $z(0) = z_0$ satisfies $|z(t, z_0)| \leq \exp(-\sigma t) a(|z_0|)$ for all $t \geq 0$. Suppose furthermore that $0 \in \Re^n$ is robustly K-exponentially stable for (2.2), i.e., there exist constants $\sigma > 0$, $\lambda \in (0,1)$ and a function $a \in K_\infty$ such that for every $x_0 \in \Re^n$ and for every locally bounded input $u : \Re^+ \to [0,+\infty)$ the solution $x(t, x_0; u)$ of (2.2) with initial condition $x(0) = x_0$ corresponding to $u : \Re^+ \to [0,+\infty)$ satisfies $|x(t, x_0; u)| \leq \exp(-\lambda \sigma t) a(|x_0|)$ for all $t \geq 0$.

The global discretization error $e(t, x_0; u)$ at time $t \geq 0$ for the numerical solution with initial condition $x(0) = x_0$ corresponding to $u : \Re^+ \to [0,+\infty)$ is the difference $e(t, x_0; u) := z(t, x_0) - x(t, x_0; u)$. Clearly we have $|e(t, x_0; u)| \leq 2 \exp(-\lambda \sigma t) a(|x_0|)$ for all $t \geq 0$, i.e., solvability of problem (P1) implies that the global discretization error $e(t, x_0; u)$ tends to 0 as $t \to +\infty$.

Define $d(s, h, x_0) := z(s, x_0) - x_0 - sF(h, x_0)$ (the local discretization error). We have

$$\begin{aligned} e(\tau_i + s, x_0; u) &= z(\tau_i + s, x_0) - x(\tau_i + s, x_0; u) \\ &= z(s, z(\tau_i, x_0)) - z(\tau_i, x_0) + z(\tau_i, x_0) - sF(h_i, x(\tau_i, x_0; u)) - x(\tau_i, x_0; u) = \\ &= d(s, h_i, z(\tau_i, x_0)) + sF(h_i, z(\tau_i, x_0)) - sF(h_i, x(\tau_i, x_0; u)) + e(\tau_i, x_0; u) \end{aligned} \quad (7.1)$$

Suppose that there exist a continuous function $L: \Re^n \to \Re^+$ such that

$$|F(h, z) - F(h, x)| \leq L(x_0)|z - x|,$$
for all $z, x \in \Re^n$ with $|x| \leq a(|x_0|)$, $|z| \leq a(|x_0|)$ and $h \in [0, \varphi(x)]$ \quad (7.2)

Notice that for implicit Runge-Kutta methods inequality (7.2) may not be satisfied since $F(h, z)$ may not be defined for all $h \in [0, \varphi(x)]$. However, for explicit Runge-Kutta methods inequality (7.2) holds in general for appropriate continuous function $L: \Re^n \to \Re^+$.

Define $\tilde{d}(h, x) := \dfrac{d(h, h, x)}{h} = \dfrac{z(h, x) - x}{h} - F(h, x)$. Standard arguments show that the following inequality holds for all $i \geq 0$:

$$|e(\tau_{i+1}, x_0; u)| \leq \exp(L(x_0)\tau_{i+1}) \sum_{j=0}^{i} \exp(-L(x_0)\tau_{j+1}) h_j |\tilde{d}(h_j, z(\tau_j, x_0))| \quad (7.3)$$

Noticing that $\exp(-L(x_0)\tau_{j+1}) h_j \leq \exp(-L(x_0)\tau_{j+1}) \int_{\tau_j}^{\tau_{j+1}} ds \leq \int_{\tau_j}^{\tau_{j+1}} \exp(-L(x_0)s) ds$, (7.3) implies the following inequality for all $i \geq 0$:

$$|e(\tau_{i+1}, x_0; u)| \leq \dfrac{D_i(x_0)}{L(x_0)} (\exp(L(x_0)\tau_{i+1}) - 1) \quad (7.4)$$

where

$$D_i(x_0) := \max_{j=0,\ldots,i} |\tilde{d}(h_j, z(\tau_j, x_0))| \quad (7.5)$$

and therefore

$$|e(\tau_{i+1}, x_0; u)| \leq \min\left\{ \dfrac{D_i(x_0)}{L(x_0)} (\exp(L(x_0)\tau_{i+1}) - 1), 2\exp(-\lambda \sigma \tau_{i+1}) a(|x_0|) \right\}, \quad \forall i \geq 0$$

The above inequality leads to the following estimate of the global discretization error:

$$|e(\tau_{i+1}, x_0; u)| \leq \left( \dfrac{D_i(x_0)}{L(x_0)} \right)^{\frac{\lambda \sigma}{\lambda \sigma + L(x_0)}} (2 a(|x_0|))^{\frac{L(x_0)}{\lambda \sigma + L(x_0)}}, \quad \forall t \geq 0 \quad (7.6)$$

Inequality (7.6) has important consequences:

**a)** Inequality (7.6) shows an order reduction phenomenon. For example, if the Runge-Kutta scheme used is of order $p$ then $D_i(x_0) = O(\overline{h}_i^p)$, where $\overline{h}_i := \max_{j=0,\ldots,i} h_j$, and consequently we obtain from (7.6)

$$\sup_{i \geq 0} |e(\tau_i, x_0; u)| \leq O(\overline{h}^{\frac{p \lambda \sigma}{\lambda \sigma + L(x_0)}}), \text{ where } \overline{h} := \max_{j \geq 0} h_j.$$

**b)** A second important point that should be emphasized is the possibility of controlling the global discretization error. The demand $|e(\tau_i, x_0; u)| \leq \varepsilon$ for all $i \geq 1$, where $\varepsilon > 0$ is given, is guaranteed by (7.6) if

$$D_i(x_0) \leq L(x_0) \varepsilon^{1 + \frac{q(x_0)}{\lambda}} (2a(|x_0|))^{-\frac{q(x_0)}{\lambda}}, \text{ for all } i \geq 0$$

where $q(x_0) := \frac{L(x_0)}{\sigma}$, or equivalently

$$\left|\frac{z(h_i + \tau_i, x_0) - z(\tau_i, x_0)}{h_i} - F(h_i, z(\tau_i, x_0))\right| \leq L(x_0)\varepsilon^{1+\frac{q(x_0)}{\lambda}}\left(2a(|x_0|)\right)^{-\frac{q(x_0)}{\lambda}}, \text{ for all } i \geq 0$$

For a Runge-Kutta method of order $p \geq 1$, there exists a continuous function $K : \Re^n \to (0, +\infty)$ (see (2.8a,b) and use of the fact that $f$ is locally Lipschitz) such that:

$$\left|\frac{z(h_i + \tau_i, x_0) - z(\tau_i, x_0)}{h_i} - F(h_i, z(\tau_i, x_0))\right| \leq h_i^p K(x_0) \max_{0 \leq s \leq h_i} |z(\tau_i + s, x_0)| \quad (7.7)$$

Therefore, by virtue of (7.7), the demand $|e(\tau_i, x_0; u)| \leq \varepsilon$ for all $i \geq 1$ is guaranteed if

$$h_i \leq \min\left\{\left(\frac{2L(x_0)}{K(x_0)}\right)^{\frac{1}{p}} \exp\left(\frac{\sigma}{p}\tau_i\right)\left(2\varepsilon^{-1}a(|x_0|)\right)^{-\frac{q(x_0)+\lambda}{p\lambda}}, \varphi(x(\tau_i, x_0; u))\right\}, \text{ for all } i \geq 0 \quad (7.8)$$

For example, using the explicit Euler method we get $p = 1$, $K(x_0) := \frac{1}{2}L^2(x_0)$. Thus we obtain from (7.8):

$$h_i \leq \min\left\{\frac{4}{L(x_0)}\exp(\sigma \tau_i)\left(2\varepsilon^{-1}a(|x_0|)\right)^{-\frac{q(x_0)+\lambda}{\lambda}}, \varphi(x(\tau_i, x_0; u))\right\}$$

Similar computations can be made for all Runge-Kutta schemes. Notice that as time becomes larger (and consequently the solution approaches the equilibrium point), larger time steps are allowed.

## 8. Conclusions

In this work, we considered the problem of step size selection for numerical schemes such that the numerical solution presents the same qualitative behavior as the original system of ODEs. Specifically, we developed tools for nonlinear systems with a globally asymptotically stable equilibrium point which are similar to methods used in Nonlinear Control Theory. It is shown how the problem of appropriate step size selection can be converted to a rigorous abstract feedback stabilization problem for a particular hybrid system. Feedback stabilization methods based on Lyapunov functions and Small-Gain results were employed. The obtained results have been applied to several examples and have been shown to be efficient for controlling the global discretization error.

The methodology presented in the present work (transformation of the step selection problem to a feedback stabilization problem and use of modern nonlinear control theory for the solution of the problem) can be used for more complicated numerical problems such as the step size selection problem for:

- ➢ the numerical approximation of the solution of infinite-dimensional systems (systems described by partial differential equations or systems described by retarded functional differential equations),
- ➢ systems with more complicated attractors,
- ➢ time-varying systems,
- ➢ systems with inputs.

Future work will address these problems.

# Appendix

**Proof of Lemma 4.3:**

(i) Define $Q(x) := -\nabla V(x) f(x)$. The function $Q : \Re^n \to \Re^+$ is continuous and by virtue of (4.7) is positive definite too. Standard results on inf-convolutions guarantee that the function $W(x) := \inf\{Q(y) + |y - x| : y \in \Re^n\}$ is globally Lipschitz on $\Re^n$ (with unit Lipschitz constant), positive definite and satisfies (4.8).

(ii) Since $f(0) = 0$, it follows that for all $x \in \Re^n$

$$|f(z)| \le l_f(x)|z|, \text{ for all } z \in \Re^n \text{ with } V(z) \le V(x) \tag{A1}$$

Inequality (A1) in conjunction with the fact that $V(z(t,x)) \le V(x)$ for all $t \ge 0$ and Gronwall's inequality implies

$$\exp(-l_f(x)t)|x| \le |z(t,x)| \le \exp(l_f(x)t)|x|, \text{ for all } (t,x) \in \Re^+ \times \Re^n \tag{A2}$$

Therefore definition (4.10) and inequality (A2) imply that

$$\exp(-b)|x| \le |z(h,x)| \le \exp(b)|x|, \text{ for all } h \in [0, \varphi(x)] \tag{A3}$$

Let $W : \Re^n \to \Re^+$ be the locally Lipschitz, positive definite function which satisfies inequality (4.8). Define:

$$\widetilde{W}(x) := \min\{W(y) : y \in \Re^n, \exp(-b)|x| \le |y| \le \exp(b)|x|\} \tag{A4}$$

Clearly, definition (A4) guarantees that $\widetilde{W} : \Re^n \to \Re^+$ is a continuous, positive definite function. Moreover, by virtue of inequalities (4.8), (A3) and definition (A4) we obtain for all $h \in [0, \varphi(x)]$ and $x \in \Re^n$:

$$V(z(h,x)) - V(x) = \int_0^h \nabla V(z(s,x)) f(z(s,x)) ds \le -\int_0^h W(z(s,x)) ds \le -h \widetilde{W}(x) \tag{A5}$$

i.e., the desired inequality (4.9).

(iii) Notice first that inequality (4.11) guarantees that there exists $\varphi \in C^0(\Re^n; (0, +\infty))$ satisfying (4.12). Define:

$$M_f^b(x) := \max\{|f(y)| : y \in \Re^n, |y| \le \exp(b)|x|\} \tag{A6}$$

Notice that inequality (A1) and definition (A6) imply that

$$M_f^b(x) \le l_f(x) \exp(b)|x| \tag{A7}$$

Taking into account inequalities (A3), (A7), (4.12) in conjunction with definition (A6) we obtain for all $h \in [0, \varphi(x)]$ and $x \in \Re^n$:

$$|W(x) - W(z(h,x))| \le l_W^b(x)|x - z(h,x)| \le l_W^b(x) l_f(x) \exp(b)|x| \varphi(x) \tag{A8}$$

Inequalities (A8) and (4.12) imply the following inequality for all $h \in [0, \varphi(x)]$ and $x \in \Re^n$:

$$-W(z(h,x)) \le -\lambda W(x) \tag{A9}$$

Moreover, by virtue of inequalities (4.8), (A3) and definition (A4) we obtain for all $h \in [0, \varphi(x)]$ and $x \in \Re^n$:

$$V(z(h,x)) - V(x) = \int_0^h \nabla V(z(s,x)) f(z(s,x)) ds \leq -\int_0^h W(z(s,x)) ds \leq -\lambda h W(x) \tag{A10}$$

i.e., the assertion. The proof is complete. ◁

**Proof of Proposition 4.4:** Since $0 \in \Re^n$ is locally exponentially stable for (2.1), it follows that the matrix $A := Df(0)$ is Hurwitz. Consequently there exists a symmetric, positive definite matrix $P \in \Re^{n \times n}$ and a constant $\mu > 0$ such that

$$x'(A'P + PA)x \leq -2\mu |x|^2, \quad \forall x \in \Re^n \tag{A11}$$

(see ...). Consequently, there exists constant $\delta > 0$ (sufficiently small) such that

$$2x'Pf(x) \leq -\mu |x|^2, \text{ for all } x \in \Re^n \text{ with } |x| \leq 2\delta \tag{A12}$$

Let $W : \Re^n \to \Re^+$ a continuously differentiable function with

$$W(x) := \begin{cases} -2x'Pf(x) & \text{for } |x| < \delta \\ \mu |x|^2 \left(1 + |f(x)|^2\right) & \text{for } |x| > 2\delta \end{cases} \text{ and } W(x) \geq \mu |x|^2 \text{ for all } x \in \Re^n \tag{A13}$$

and define the function

$$V(x) := \int_0^{+\infty} W(z(t,x)) dt \tag{A14}$$

By virtue of Theorem 2.46 in [4] $V$ as defined by (A14) is a Lyapunov function for (2.1) satisfying

$$\nabla V(x) f(x) = -W(x) \text{ for all } x \in \Re^n \tag{A15}$$

By virtue of Proposition 2.48 in [4] is the unique function satisfying (A15) with $V(0) = 0$. An inspection of the proof yields that if equation (A15) holds on a forward invariant set for (2.1) then uniqueness holds on this set (because uniqueness is established by looking at trajectories in forward time). Notice that by virtue of (A15) and definition (A13), it follows that (4.14) holds.

Picking a forward invariant neighborhood $N \subset B_\delta(0)$ of zero (such a forward invariant neighborhood of zero exists because $0 \in \Re^n$ is asymptotically stable), we observe by (A13) that the function $\tilde{V}(x) = x'Px$ satisfies (A15) as well on $N \subset B_\delta(0)$ and $\tilde{V}(0) = 0$. Consequently, $V(x) \equiv \tilde{V}(x) = x'Px$ on $N \subset B_\delta(0)$. Let $\varepsilon > 0$ sufficiently small with $B_\varepsilon(0) \subseteq N$. It follows that (4.13) holds. The proof is complete. ◁